\documentclass[11pt]{article}
\usepackage{amssymb,amsthm,hyperref}

\newtheorem{thm}{Theorem}[section]
 
 \newtheorem{cor}{Corollary}[section]
 \newtheorem{lem}{Lemma}[section]
 \newtheorem{prop}{Proposition}[section]
 \newtheorem{defn}{Definition}[section]
\theoremstyle{remark}
\newtheorem{rem}{Remark}[section]

\title{\textbf{Existence and asymptotic behavior of $C^1$ solutions \\ to the multidimensional compressible \\ Euler equations with damping} \thanks {This work is supported by
NSFC 10571158. }}
\author{\textbf{Daoyuan Fang}\thanks {E-mail: dyf@zju.edu.cn},\ \ \ \ \textbf{Jiang Xu}\thanks {E-mail: jiangxu\underline{ }79@yahoo.com.cn}\\
\small{\textsl{Department of Mathematics, Zhejiang University},}
\\ \small{\textsl{Hangzhou 310027, P.R.China}}}
\date{}
\begin{document}
\maketitle{} \hspace{-5mm}\textbf{Abstract}\\

 \small{In this paper, the existence and asymptotic behavior of $C^1$ solutions to the multidimensional compressible Euler equations with damping
   on the framework of Besov space are considered. We \textit{weaken} the regularity requirement of the initial data, and improve the well-posedness results of \textsc{Sideris}-\textsc{Thomases}-\textsc{Wang} (Comm.P.D.E. 28 (2003) 953).
   The global existence
 lies on a crucial a-priori estimate which is proved by the \textit{spectral localization} method. The main analytic tools are the
Littlewood-Paley decomposition and Bony's para-product
 formula}.\\

\hspace{-0.5cm}\textsl{MSC:} \small{35L65;\ 76N15}\\

\hspace{-0.5cm}\textsl{Keywords}: \small{Euler equations; damping;
classical
 solutions; spectral localization}

\section{Introduction and main results}
\hspace{5mm}In this paper, we study the following Euler equation
with damping for a perfect gas flow:
$$\cases{n_{t}+\nabla\cdot(n\textbf{u})=0\cr
 (n\textbf{u})_{t}+\nabla\cdot(n\textbf{u}\otimes\textbf{u})+\nabla p(n)=-an\textbf{u} }\eqno(1.1)$$
 for $ (t,x)\in[0,+\infty)\times\mathbb{R}^{N}, \ N\geqslant1$,
 where $n$ and $\textbf{u}=(u^1,u^2,\cdot\cdot\cdot,u^{N})^{T}$($T$ represents transpose) denote the density, velocity for
 the gas respectively. $n\textbf{u}$ stands for the momentum. The pressure
 $p$ satisfies the $\gamma$-law: $$p=p(n)=An^{\gamma},\eqno(1.2)$$ where the case $\gamma>1$  corresponds to the isentropic gas and $\gamma=1$ corresponds to the isothermal gas, $A$ is
 a positive
 constant. The positive constant $a$ is the
 damping coefficient. The system
 is supplemented with the initial data
 $$(n,\textbf{u})(x,0)=(n_{0},\textbf{u}_{0})(x),\ \
  x\in\mathbb{R}^{N}. \eqno(1.3)$$

The system (1.1) describes that the compressible gas flow passes a
porous medium and the medium induces a friction force, proportional
to the linear momentum in the opposite direction. It is hyperbolic
with two characteristic speeds $\lambda=\textbf{u}\pm\sqrt{p'(n)}$.
As a vacuum appears, it fails to be strict hyperbolic. Thus, the
system involves three mechanisms: nonlinear convection, lower-order
dissipation of damping and the resonance due to vacuum. After
\textsc{Nishida}'s \cite{N1,N2} pioneer works for (1.1), many
contributions have been made on the small smooth solutions and
piecewise smooth Riemann solutions away from vacuum, we can cite
\cite{HL,H,HP1,HS,WY,STW} and their references. Among them, for the
one dimensional case \cite{HL}, the system can be written in the
Lagrangian coordinates as follows:
$$\cases{v_{t}-u_{x}=0,\cr
 u_{t}+ p(v)_{x}=-au, }\eqno(1.4)$$
where $v=1/n$ is the specific volume. It was shown that the system
(1.4) was time asymptotically equivalent to the porous media
equation. In \cite{STW}, \textsc{Sideris}, \textsc{Thomases }\&
\textsc{Wang} explained that the damping only presented weak
dissipation in 3D space: it could prevent the development of
singularities if the initial data was small and smooth,
furthermore, they obtained the decay of classical solutions to the
constant background state in $L^{\infty}$ at a rate of $
(1+t)^{-(3/2)}$, but singularities was exhibited for large data
under some assumption. However, the main open problems for (1.1)
with vacuum are still far from well-known. One of them is to study
the singular evolution of the vacuum interface. As the first step
in this direction, \textsc{Xu} \& \textsc{Yang} \cite{XY} proved
that a local existence theorem on a perturbation of a planar wave
solution for (1.1) under the assumption of physical vacuum
boundary condition. For the large-time asymptotic behavior for the
solutions with vacuum, recently, \textsc{Huang} \& \textsc{Pan}
\textit{et al}. \cite{HMP,HP2}  gave a complete answer to this
problem. In fact, they showed that the $L^{\infty}$ weak entropy
solutions with vacuum for the Cauchy problem converged to the
Barenblatt's profile of the porous medium equation strongly in
$L^{p}$.

In this paper, we are concerned with lowering the regularity of the
initial data in the generally multidimensional space. As in
\cite{STW}, we also consider a perturbation of the constant
equilibrium state $(\bar{n},0)(\bar{n}>0)$. First of all, we give a
local existence result in Besov space
$B^{\sigma}_{2,1}(\mathbb{R}^{N})$ ($\sigma=1+\frac{N}{2}$) for
(1.1)-(1.3) away from vacuum.
\begin{thm}($N\geqslant1$)
Suppose that\ $(n_{0}-\bar{n}, \textbf{u}_{0})\in
B^{\sigma}_{2,1}(\mathbb{R}^{N})$ with \ $n_{0}>0$, then there
exist a time $T_{0}>0$ and a unique solution $(n,\textbf{u})$ of
the system (1.1)-(1.3) such that
$$(n,\textbf{u})\in C^{1}([0,T_{0}]\times \mathbb{R}^{N})\ \ \ \mbox{with}\
\ \  n>0\ \  \mbox{for all} \ \ t\in [0,T_{0}]$$ and
$$(n-\bar{n},\textbf{u})\in  C([0,T_{0}],B^{\sigma}_{2,1}(\mathbb{R}^{N}))\cap
C^1([0,T_{0}],B^{\sigma-1}_{2,1}(\mathbb{R}^{N})).$$
\end{thm}
\begin{rem}
The nonlinear pressure term makes computation more fussy in the
\textit{spectral localization} estimates due to commutors. To get
around this, we introduce a new variable (sound speed) which
transforms the nonlinear term into linear and double-linear terms in
virtue of the ideas in \cite{STW}. In fact, the original system
(1.1)-(1.3) transforms into a symmetric hyperbolic system
(3.1)-(3.2) where we can obtain the effective a-priori estimates.
Different from the local existence result in \cite{STW}, Theorem 1.1
follows from Proposition 4.1, Remark 4.1 and Remark 3.1. The proof
of Proposition 4.1 is organized as follows. First, we regularize the
initial data of (3.1)-(3.2) and obtain approximative local solutions
based on Kato's results. Second, we find a uniform positive time
$T_{0}$ such that the approximative solution sequence is uniform
bounded in $C([0,T_{0}];B^{\sigma}_{2,1}(\mathbb{R}^{N}))\cap
C^1([0,T_{0}],B^{\sigma-1}_{2,1}(\mathbb{R}^{N}))$. Finally,  we
utilize the compactness argument to pass the limit (For detail, see
Proposition 4.1.).
\end{rem}
Under a smallness assumption, we establish the global existence of
classical solutions in Besov space
$B^{\sigma+\varepsilon}_{2,2}(\mathbb{R}^{N})$
($\sigma=1+\frac{N}{2},\ \varepsilon>0$) for (1.1)-(1.3).
\begin{thm}($N\geqslant3$) Suppose that \ $(n_{0}-\bar{n}, \textbf{u}_{0})\in
B^{\sigma+\varepsilon}_{2,2}(\mathbb{R}^{N})$. There exists a
positive constant $\delta_{0}$ depending only on $A, \gamma, a$
and $\bar{n}$ such that if
$$\|(n-\bar{n},\textbf{u})(\cdot,0)\|^2_{B^{\sigma+\varepsilon}_{2,2}(\mathbb{R}^{N})}+\|(n_{t},\textbf{u}_{t})(\cdot,0)\|^2_{B^{\sigma-1+\varepsilon}_{2,2}(\mathbb{R}^{N})}\leqslant
\delta_{0},\eqno(1.5)$$ then there exists a unique global solution
$(n,\textbf{u})$ of the system (1.1)-(1.3) satisfying
$$(n,\textbf{u})\in C^{1}([0,\infty)\times \mathbb{R}^{N})$$
and
$$(n-\bar{n},\textbf{u}) \in
C([0,\infty),B^{\sigma+\varepsilon}_{2,2}(\mathbb{R}^{N}))\cap
C^1([0,\infty),B^{\sigma-1+\varepsilon}_{2,2}(\mathbb{R}^{N})).$$
Moreover, we have the energy estimate\begin{eqnarray*}
&&\|(n-\bar{n},\textbf{u})(\cdot,t)\|^2_{B^{\sigma+\varepsilon}_{2,2}(\mathbb{R}^{N})}+\|(n_{t},\textbf{u}_{t})(\cdot,t)\|^2_{B^{\sigma-1+\varepsilon}_{2,2}(\mathbb{R}^{N})}\\&&+\mu_{0}\int_{0}^{t}(\|\textbf{u}(\cdot,\tau)\|^2_{B^{\sigma+\varepsilon}_{2,2}(\mathbb{R}^{N})}
+\|(\nabla
n,n_{t},\textbf{u}_{t})(\cdot,\tau)\|^2_{B^{\sigma-1+\varepsilon}_{2,2}(\mathbb{R}^{N})})d\tau\\
&\leqslant&
\|(n-\bar{n},\textbf{u})(\cdot,0)\|^2_{B^{\sigma+\varepsilon}_{2,2}(\mathbb{R}^{N})}+\|(n_{t},\textbf{u}_{t})(\cdot,0)\|^2_{B^{\sigma-1+\varepsilon}_{2,2}(\mathbb{R}^{N})},\
\ \ t\geqslant0,\end{eqnarray*} where the positive constant
$\mu_{0}$ depends only on $A, \gamma, a$ and $\bar{n}$.
\end{thm}
\begin{rem} In fact, the smallness of $
(n_{t},\textbf{u}_{t})(x,0)$ can be derived by Eqs.(1.1) and the
smallness of $(n-\bar{n},\textbf{u})(x,0)$. For the simplicity of
the statement, we give the assumption (1.5) directly.
\end{rem}

\begin{rem} Theorem 1.2 follows from Proposition 5.1 and Remark 3.1. The proof of a crucial
a-priori estimate (Proposition 5.2) is separated into the low
frequency part (Lemma 5.5) and high frequency part (Lemma 5.6)
elaborately. On each high frequency ($q\geqslant0$),
$\|\triangle_{q}\nabla m\|_{L^2}$ is equivalent to
$2^{q}\|\triangle_{q} m\|_{L^2}$ by Lemma 2.1, but it is not valid
for low frequency ($q=-1$). In \cite{FXZ}, we knew that the
Poisson potential remedied the estimate on
$\|\triangle_{-1}m\|_{L^2}$. Here, we can't obtain any estimates
on $\|\triangle_{-1}m\|_{L^2}$, however, with the help of
H\"{o}lder's inequality and Gagliardo-Nirenberg-Sobolev inequality
($N>2$), we can get the estimates on $\|\triangle_{-1}\nabla
m\|_{L^2}$ (For detail, see (5.12), (5.18), (5.20), (5.24) and
(5.26)). Hence, in order to ensure our functional space still
imbedding into $C^{1}$ space, we need to increase a little of
regularity, furthermore, which leads to the global existence of
classical solutions to (3.1)-(3.2).
\end{rem}
In \cite{STW}, the authors obtained the decay of classical solutions
to the equilibrium state $(\bar{n},0)$ in $L^{\infty}$ at a rate of
$ (1+t)^{-(3/2)}$. According to the energy estimate in Theorem 1.2,
we also see the large-time asymptotic behavior of solutions in Besov
space roughly.
\begin{cor}($N\geqslant3$)
Let $(n, \textbf{u})$ be the solution in Theorem 1.2, we have
($\sigma=1+\frac{N}{2},\ \varepsilon'<\varepsilon$.)
$$\|n(\cdot,t)-\bar{n}\|_{B^{\sigma-1+\varepsilon'}_{p,2}(\mathbb{R}^{N})}\rightarrow 0 \ \ (p=\frac{2N}{N-2}),\ \ \ \ \|\textbf{u}(\cdot,t)\|_{B^{\sigma+\varepsilon'}_{2,2}(\mathbb{R}^{N})}\rightarrow 0,\ \ \mbox{as}\ \ t\rightarrow +\infty.   $$
\end{cor}
Different from the result in \cite{STW}, the following theorem
characterizes the exponential decay of the vorticity in Besov space
$B^{\sigma-1}_{2,1}(\mathbb{R}^{N})$.
\begin{thm} (N=3) Let $(n, \textbf{u})$ be the solution in Theorem 1.2. If
$$\|(n-\bar{n},\textbf{u})(\cdot,0)\|^2_{B^{\sigma+\varepsilon}_{2,2}(\mathbb{R}^{N})}+\|(n_{t},\textbf{u}_{t})(\cdot,0)\|^2_{B^{\sigma-1+\varepsilon}_{2,2}(\mathbb{R}^{N})}\leqslant
\delta_{0}',$$ then the vorticity $\omega=\nabla\times\textbf{u}$
decays exponentially in $B^{\sigma-1}_{2,1}(\mathbb{R}^{N})$:
$$\|\omega(\cdot,t)\|_{B^{\sigma-1}_{2,1}(\mathbb{R}^{N})}\leqslant \|\omega(\cdot,0)\|_{B^{\sigma-1}_{2,1}(\mathbb{R}^{N})}\exp(-\mu'_{0}t), \ \ t\geqslant 0,$$
where \ $\omega(x,0)=\nabla\times \textbf{u}_{0}$,\ the positive
constants \ $\delta_{0}'=\min\{\delta_{0}, \frac{a^2}{4C^2_{0}}\}$
and $\mu_{0}'$ depend only on $A, \gamma, a$ and $\bar{n}$ \
($C_{0}$  a uniform constant given in (5.33)).
\end{thm}

The paper is arranged as follows. In Section 2, we present some
definitions and basic facts on the Littlewood-Paley decomposition
and Bony's para-product formula. In Section 3, we reformulate the
system (1.1)-(1.3) in order to obtain the effective a-priori
estimates by the spectral localization method. In Section 4, we
are concerned with the local existence and uniqueness of classical
solutions to (3.1)-(3.2) with general initial data. In the last
section, we deduce a crucial a-priori estimate under a smallness
assumption which is used to complete the proof of global
existence. Finally, it is also shown that the vorticity decays to
zero in time exponentially.

Throughout this paper, the symbol $C$ denotes a harmless constant.
All functional spaces will be considered in $\mathbb{R}^{N}$, so we
can omit the space dependence for simplicity.

\section{Littlewood-Paley analysis}
In this section, these definitions and basic facts can be found in
\textsc{Darchin's} mini-course \cite{D}.

 Let $\mathcal{S}(\mathbb{R}^{N})$ be the Schwarz class. ($\varphi, \chi)$ is a couple of smooth functions valued in [0,1]
such that $\varphi$ is supported in the shell
$C(0,\frac{3}{4},\frac{8}{3})=\{\xi\in\mathbb{R}^{N}|\frac{3}{4}\leqslant|\xi|\leqslant\frac{8}{3}\}$,
$\chi$ is supported in the ball $B(0,\frac{4}{3})=
\{\xi\in\mathbb{R}^{N}||\xi|\leqslant\frac{4}{3}\}$ and $$
\chi(\xi)+\sum_{q=0}^{\infty}\varphi(2^{-q}\xi)=1,\ \ \ \ q\in
\mathbb{Z},\ \  \xi\in\mathbb{R}^{N}.$$ For
$f\in\mathcal{S'}$(denote the set of temperate distributes which
is the dual one of $\mathcal{S}$), we can define the
nonhomogeneous dyadic blocks as follows:
$$\Delta_{-1}f:=\chi(D)f=\tilde{h}\ast f\ \ \ \mbox{with}\ \ \tilde{h}=\mathcal{F}^{-1}\chi,$$
$$\Delta_{q}f:=\varphi(2^{-q}D)f=2^{qN}\int h(2^{q}y)f(x-y)dy\ \ \ \mbox{with}\ \ h=\mathcal{F}^{-1}\varphi,\ \ \mbox{if}\ \ q\geqslant0,$$
where $\ast, \ \ \mathcal{F}^{-1} $ represent the convolution
operator and the inverse Fourier transform, respectively. The
nonhomogeneous Littlewood-Paley decomposition is$$ f=\sum_{q
\geqslant-1}\Delta_{q}f \ \ \ \mbox{in}\ \ \ \mathcal{S'}.$$
Define the low frequency cut-off by
$$S_{q}f:=\sum_{p\leqslant q-1}\Delta_{p}f.$$ Of course, $S_{0}f=\triangle_{-1}f$. The above
Littlewood-Paley decomposition is almost orthogonal in $L^2$.
\begin{prop}
For any $f\in\mathcal{S'}(\mathbb{R}^{N})$ and
$g\in\mathcal{S'}(\mathbb{R}^{N})$, the following properties hold:
$$\Delta_{p}\Delta_{q}f\equiv 0 \ \ \ \mbox{if}\ \ \ |p-q|\geqslant 2,$$
$$\Delta_{q}(S_{p-1}f\Delta_{p}g)\equiv 0\ \ \ \mbox{if}\ \ \ |p-q|\geqslant 5.$$
\end{prop}
Besov space can be characterized in virtue of the Littlewood-Paley
decomposition.
\begin{defn}
Let $1\leqslant p\leqslant\infty$ and $s\in \mathbb{R}$. For $1\leqslant r<\infty$, the Besov spaces  $B^{s}_{p,r}(\mathbb{R}^{N})$ are defined by
$$f\in B^{s}_{p,r}(\mathbb{R}^{N}) \Leftrightarrow \Big(\sum_{q\geqslant-1}(2^{qs}\|\Delta_{q}f\|_{L^{p}})^{r}\Big)^{\frac{1}{r}}<\infty$$
and $B^{s}_{p,\infty}(\mathbb{R}^{N})$ are defined by
$$f\in B^{s}_{p,\infty}(\mathbb{R}^{N}) \Leftrightarrow \sup_{q\geqslant-1}2^{qs}\|\Delta_{q}f\|_{L^{p}}<\infty.$$
In particular, $B^{s}_{2,2}(\mathbb{R}^{N})\equiv
H^{s}(\mathbb{R}^{N})$.
\end{defn}
\begin{defn}
Let $f,g $ be two temperate distributions. The product $f\cdot g$ has the Bony's decomposition:
$$f\cdot g=T_{f}g+T_{g}f+R(f,g), $$
where $T_{f}g$ is paraproduct of $g$ by $f$,
$$ T_{f}g=\sum_{p\leqslant q-2}\Delta_{p}f\Delta_{q}g=\sum_{q}S_{q-1}f\Delta_{q}v$$
and the remainder $ R(f,g)$ is denoted by
$$R(f,g)=\sum_{q}\Delta_{q}f\tilde{\Delta}_{q}g\ \ \ \mbox{with} \ \
\tilde{\Delta}_{q}:=\Delta_{q-1}+\Delta_{q}+\Delta_{q+1}.$$
\end{defn}
\begin{lem}(Bernstein)
Let $k\in\mathbb{N}$ and $0<R_{1}<R_{2}$. There exists a constant
$C$ depending only on $R_{1},R_{2}$ and $N$ such that for all
$1\leqslant a\leqslant b\leqslant\infty$ and $f\in L^{a}$, we have
$$\mathrm{Supp}\ \mathcal{ F}f\subset \textsl{B}(0,R_{1}\lambda)\Rightarrow\sup_{|\alpha|=k}\|\partial^{\alpha}f\|_{L^{b}}\leqslant C^{k+1}\lambda^{k+N(\frac{1}{a}-\frac{1}{b})}\|f\|_{L^{a}};$$
$$\mathrm{Supp}\ \mathcal{ F}f\subset \textsl{C}(0,R_{1}\lambda,R_{2}\lambda)\Rightarrow C^{-k-1}\lambda^{k}\|f\|_{L^{a}}\leqslant \sup_{|\alpha|=k}\|\partial^{\alpha}f\|_{L^{a}}\leqslant C^{k+1}\lambda^{k}\|f\|_{L^{a}}.$$
Here, $\mathcal{F}$ represents the Fourier transform.
\end{lem}
A result of compactness for Besov space:
\begin{prop}
Let $1\leqslant p,r\leqslant \infty,\ s\in \mathbb{R}$ and
$\varepsilon>0$. For all $\phi\in C_{c}^{\infty}$, the map
$f\mapsto\phi f$ is compact from
$B^{s+\varepsilon}_{p,r}(\mathbb{R}^{N})$ to
$B^{s}_{p,r}(\mathbb{R}^{N})$.
\end{prop}
Finally, we state a result of continuity for the composition which
is used to end this section.
\begin{prop}
Let $1\leqslant p,r\leqslant \infty$,\ \ $I$ be open interval of
$\mathbb{R}$. Let $s>0$ and $n$ be the smallest integer such that
$n\geqslant s$. Let $F:I\rightarrow\mathbb{R}$ satisfy $F(0)=0$
and $F'\in W^{n,\infty}(I;\mathbb{R}).$ Assume that $v\in
B^{s}_{p,r}$ takes values in $J\subset\subset I$. Then $F(v)\in
B^{s}_{p,r}$ and there exists a constant $C$ depending only on
$s,I,J$ and $N$ such that $$\|F(v)\|_{B^{s}_{p,r}}\leqslant
C(1+\|v\|_{L^{\infty}})^{n}\|F'\|_{W^{n,\infty}(I)}\|v\|_{B^{s}_{p,r}}.$$
\end{prop}
\section{Reformulation of the original system}
In this section, we are going to reformulate (1.1)-(1.3) in order
to get the effective a-priori estimates by the spectral localization
method. For the isentropic case $(\gamma>1)$:
 introducing the sound speed
$$\psi(n)=\sqrt{p'(n)},$$ and denoting the sound speed at a background density $\bar{n}$ by  $\bar{\psi}=\psi(\bar{n})$, as in \cite{STW}, we define
$$m=\frac{2}{\gamma-1}(\psi(n)-\bar{\psi}).$$
Then the system (1.1) is transformed into the following system for
$C^1$ solutions:
$$\cases{m_{t}+\bar{\psi}\mbox{div}\textbf{u}=-\textbf{u}\cdot\nabla m-\frac{\gamma-1}{2}m\mbox{div}\textbf{u},\cr
 \textbf{u}_{t}+\bar{\psi}\nabla m+a\textbf{u}=-\textbf{u}\cdot\nabla\textbf{u}-\frac{\gamma-1}{2}m\nabla m. }\eqno(3.1)$$
 The initial data (1.3) becomes
$$(m,\textbf{u})|_{t=0}=(m_{0},\textbf{u}_{0})\eqno(3.2)$$
with $$m_{0}=\frac{2}{\gamma-1}(\psi(n_{0})-\bar{\psi}).$$

\begin{rem} For any $T>0$,
$(n,\textbf{u})\in C^1([0,T]\times \mathbb{R}^{N})$
solves the system (1.1)-(1.2) with $n>0$, then
$(m,\textbf{u})\in C^1([0,T]\times \mathbb{R}^{N})$
solves the system (3.1)-(3.2) with
$\frac{\gamma-1}{2}m+\bar{\psi}>0$; Conversely, if
$(m,\textbf{u})\in C^1([0,T]\times \mathbb{R}^{N})$
solves the system (3.1)-(3.2) with
$\frac{\gamma-1}{2}m+\bar{\psi}>0$, let
$n=\psi^{-1}(\frac{\gamma-1}{2}m+\bar{\psi})$, then
$(n,\textbf{u})\in C^1([0,T]\times \mathbb{R}^{N})$
solves the system (1.1)-(1.2) with $n>0$.
\end{rem}
 For the isothermal case $(\gamma=1)$: set $\tilde{n}=\sqrt{A}(\ln n-\ln
\bar{n})$, then the system (1.1) can be
transformed into the following system for $C^1$ solutions:
$$\cases{\tilde{n}_{t}+\sqrt{A}\mbox{div}\textbf{u}=-\textbf{u}\cdot\nabla \tilde{n},\cr
 \textbf{u}_{t}+\sqrt{A}\nabla \tilde{n}+a\textbf{u}=-\textbf{u}\cdot\nabla\textbf{u}. }\eqno(3.3)$$
The initial data (1.3) becomes
$$(\tilde{n},\textbf{u})|_{t=0}=(\sqrt{A}(\ln n_{0}-\ln
\bar{n}), \textbf{u}_{0}). \eqno(3.4) $$
\begin{rem} For any $T>0$, $(n,\textbf{u})\in
C^1([0,T]\times \mathbb{R}^{N})$ \ solves the system (1.1)-(1.2)
with  $n>0$, \ then $(\tilde{n},\textbf{u})\in
C^1([0,T]\times \mathbb{R}^{N})$ solves the system (3.3)-(3.4);
Conversely, if $(\tilde{n},\textbf{u})\in
C^1([0,T]\times \mathbb{R}^{N})$ solves the system (3.3)-(3.4),
let $n=\bar{n}\exp(A^{-\frac{1}{2}}\tilde{n})$, then
$(n,\textbf{u})\in C^1([0,T]\times \mathbb{R}^{N})$
solves the system (1.1)-(1.2) with $n>0$.
\end{rem}
In the subsequent sections, we study the system (3.1)-(3.2) and
prove the main results in this paper only, (3.3)-(3.4) can be
studied through the similar process.
\section{Local existence}
In this section, we shall first give the estimates of some
commutors in virtue of Bony's para-product formula and the
Littlewood-Paley decomposition in Besov $B^{\sigma}_{2,1} $ space.
Second, using the regularized means and compactness argument, we
complete the proof of local existence for (3.1)-(3.2).

Applying the operator $\triangle_{q}$ to (3.1) yields
$$\cases{\partial_{t}\triangle_{q}m+(\textbf{u}\cdot\nabla)\triangle_{q}m=-\bar{\psi}\triangle_{q}\mbox{div}\textbf{u}+[\textbf{u},\triangle_{q}]\cdot\nabla
m-\frac{\gamma-1}{2}\triangle_{q}(m\mbox{div}\textbf{u}),\cr\\
 \partial_{t}\triangle_{q}\textbf{u}+(\textbf{u}\cdot\nabla)\triangle_{q}\textbf{u}+a\triangle_{q}\textbf{u}=-\bar{\psi}\triangle_{q}(\nabla m)+[\textbf{u},\triangle_{q}]\cdot\nabla\textbf{u}-\frac{\gamma-1}{2}\triangle_{q}(m\nabla
 m),}\eqno(4.1)$$
where the commutor $[f,g]=fg-gf.$ \\  Multiplying the first equation
of Eqs.(4.1) by $\triangle_{q}m$ and the second one by
$\triangle_{q}\textbf{u}$, adding the resulting equations together
and integrating them in $\mathbb{R}^{n}$, we obtain
\begin{eqnarray*}&&\frac{1}{2}\frac{d}{dt}(\|\triangle_{q}m\|^2_{L^2}+\|\triangle_{q}\textbf{u}\|^2_{L^2})+a\|\triangle_{q}\textbf{u}\|^2_{L^2}\\
&=&\frac{1}{2}\int_{\mathbb{R}^{N}}\mbox{div}\textbf{u}(|\triangle_{q}m|^2+|\triangle_{q}\textbf{u}|^2)+\int_{\mathbb{R}^{n}}([\textbf{u},\triangle_{q}]\cdot\nabla
m\triangle_{q}m+[\textbf{u},\triangle_{q}]\cdot\nabla\textbf{u}\triangle_{q}\textbf{u})
\\&&-\frac{\gamma-1}{2}\int_{\mathbb{R}^{N}}(\triangle_{q}(m\mbox{div}\textbf{u})\triangle_{q}m+\triangle_{q}(m\nabla
m)\triangle_{q}\textbf{u}).\end{eqnarray*}$$\eqno(4.2)$$

Noticing that that the bi-linear spectral localization term, we
have
$$
\int_{\mathbb{R}^{N}}[\triangle_{q}(m\mbox{div}\textbf{u})\triangle_{q}m+\triangle_{q}(m\nabla
m)\triangle_{q}\textbf{u}]=-\int_{\mathbb{R}^{N}}\triangle_{q}m(\nabla
m\cdot\triangle_{q}\textbf{u})+\int_{\mathbb{R}^{N}}[\triangle_{q},m]\nabla
m\cdot\triangle_{q}\textbf{u}+\int_{\mathbb{R}^{N}}[\triangle_{q},m]\mbox{div}\textbf{u}\triangle_{q}m.\eqno(4.3)$$
Here, we give the following lemma to estimate these commutors in
$L^2$- norm.
\begin{lem}
The following estimates hold for any $m,\textbf{u}\in
B^{\sigma}_{2,1}:$
$$2^{q\sigma}\|[m,\triangle_{q}]\mathrm{div}\textbf{u}\|_{L^2}\leqslant
Cc_{q}\|m\|_{B^{\sigma}_{2,1}}\|\textbf{u}\|_{B^{\sigma}_{2,1}},
\eqno(4.4)$$
$$2^{q\sigma}\|[m,\triangle_{q}]\nabla m\|_{L^2}\leqslant
Cc_{q}\|\nabla m\|_{L^{\infty}}\|m\|_{B^{\sigma}_{2,1}},\eqno(4.5)$$
$$2^{q\sigma}\|[\textbf{u},\triangle_{q}]\cdot\nabla m\|_{L^2}\leqslant
Cc_{q}\|\textbf{u}\|_{B^{\sigma}_{2,1}}\|m\|_{B^{\sigma}_{2,1}},\eqno(4.6)$$
$$2^{q\sigma}\|[\textbf{u},\triangle_{q}]\cdot\nabla \textbf{u}\|_{L^2}\leqslant Cc_{q}\|\nabla
\textbf{u}\|_{L^{\infty}}\|\textbf{u}\|_{B^{\sigma}_{2,1}},\eqno(4.7)$$
where $C$ denotes a  harmless  constant, $c_{q}$ denotes a sequence
such that $\|(c_{q})\|_{ {l^{1}}}\leqslant 1.$
\end{lem}
\hspace{-5mm}\textbf{Proof.}  We are going to show (4.4) holds
only, others can be proved similarly. In order to obtain (4.4) in
Besov space $B^{\sigma}_{2,1}$, we have to split $m$ into low and
high frequencies: $m=\triangle_{-1}m+\tilde{m}$. Since there
exists a radius  $0<R<\frac{3}{4}$ such that $\mathrm{Supp}\
\mathcal{F}\tilde{m}\bigcap B(0,R)=\emptyset$, Lemma 2.1 implies
$$\|\triangle_{q}\nabla\tilde{m}\|_{L^{a}}\approx
2^{q}\|\triangle_{q}\tilde{m}\|_{L^{a}}, \ \ \ a\in[1,\infty],\ \ \
q\geqslant-1.\eqno(4.8)$$ Taking advantage of Bony's decomposition,
we have
\begin{eqnarray*}
[m,\triangle_{q}]\mbox{div}\textbf{u}&=&[\tilde{m},\triangle_{q}]\mbox{div}\textbf{u}+[\triangle_{-1}m,\triangle_{q}]\mbox{div}\textbf{u}\\
&=&\tilde{m}\triangle_{q}\mbox{div}\textbf{u}-\triangle_{q}(\tilde{m}\mbox{div}\textbf{u})+[\triangle_{-1}m,\triangle_{q}]\mbox{div}\textbf{u}
\\&=&T_{\tilde{m}}\triangle_{q}\mbox{div}\textbf{u}+T_{\triangle_{q}\mbox{div}\textbf{u}}\tilde{m}+R(\tilde{m},\triangle_{q}\mbox{div}\textbf{u})\\
&&-\triangle_{q}(T_{\tilde{m}}\mbox{div}\textbf{u}+T_{\mbox{div}\textbf{u}}\tilde{m}+R(\tilde{m},\mbox{div}\textbf{u}))+[\triangle_{-1}m,\triangle_{q}]\mbox{div}\textbf{u}.
\end{eqnarray*}
Then we can write
$[m,\triangle_{q}]\mbox{div}\textbf{u}=\sum^6_{i=1}F^{i}_{q} $ where
\begin{eqnarray*}F^1_{q}&=&T_{\tilde{m}}\triangle_{q}\partial_{j}u^{j}-\triangle_{q}T_{\tilde{m}}\partial_{j}u^{j},\ \ \ \ (\mbox{div}\textbf{u}:=\partial_{j}u^{j})\\
F^2_{q}&=&T_{\triangle_{q}\partial_{j}u^{j}}\tilde{m},\\
F^3_{q}&=&-\triangle_{q}T_{\partial_{j}u^{j}}\tilde{m},\\
F^4_{q}&=&\partial_{j}R(\tilde{m},\triangle_{q}u^{j})-\partial_{j}\triangle_{q}R(\tilde{m},u^{j}),\\
F^5_{q}&=&\triangle_{q}R(\partial_{j}\tilde{m},u^{j})-R(\partial_{j}\tilde{m},\triangle_{q}u^{j})\\
F^6_{q}&=&[\triangle_{-1}m,\triangle_{q}]\mbox{div}\textbf{u}.
\end{eqnarray*}
By Proposition 2.1, we have
\begin{eqnarray*}
F^1_{q}&=&\sum_{q'}S_{q'-1}\tilde{m}\triangle_{q'}\triangle_{q}\partial_{j}u^{j}-\triangle_{q}\sum_{q'}S_{q'-1}\tilde{m}\triangle_{q'}\partial_{j}u^{j}\\
&=&\sum_{|q-q'|\leqslant4}[S_{q'-1}\tilde{m},\triangle_{q}]\partial_{j}\triangle_{q'}u^{j}\\
&=&\sum_{|q-q'|\leqslant4}\int_{\mathbb{R}^{N}}
h(y)[S_{q'-1}\tilde{m}(x)-S_{q'-1}\tilde{m}(x-2^{-q}y)]\partial_{j}\triangle_{q'}u^{j}(x-2^{-q}y)dy.
\end{eqnarray*}
Then, applying first order Taylor's formula, Young's inequality
and (4.8), we get
\begin{eqnarray*}
2^{q\sigma}\|F^1_{q}\|_{L^2}&\leqslant&C\sum_{|q-q'|\leqslant4}\|\nabla \tilde{m}\|_{L^{\infty}}2^{(\sigma-1)(q-q')}2^{q'\sigma}\|\triangle_{q'}u^{j}\|_{L^2}\\
&\leqslant&\ Cc_{q1}\|\nabla
m\|_{L^{\infty}}\|\textbf{u}\|_{B^{\sigma}_{2,1}},\ \ \ \
c_{q1}=\sum_{|q-q'|\leqslant4}\frac{2^{q'\sigma}\|\triangle_{q'}\textbf{u}\|_{L^2}}{9\|\textbf{u}\|_{B^{\sigma}_{2,1}}}.
\end{eqnarray*}
and
\begin{eqnarray*}
2^{q\sigma}\|F^2_{q}\|_{L^2}&=& 2^{q\sigma}\Big\|\sum_{q'\geqslant q-3}S_{q'-1}\partial_{j}\triangle_{q}u^{j}\triangle_{q'}\tilde{m}\Big\|_{L^2}\\
&\leqslant& 2^{q\sigma}\sum_{q'\geqslant q-3}\|\triangle_{q'}\tilde{m}\|_{L^{\infty}}\|S_{q'-1}\partial_{j}\triangle_{q}u^{j}\|_{L^2}\\
&\leqslant&C\sum_{q'\geqslant q-3}2^{q-q'}\|\nabla m\|_{L^{\infty}}2^{q\sigma}\|\triangle_{q}\textbf{u}\|_{L^2}\\
&\leqslant& Cc_{q2}\|\nabla
m\|_{L^{\infty}}\|\textbf{u}\|_{B^{\sigma}_{2,1}},\ \ \
c_{q2}=\frac{2^{q\sigma}\|\triangle_{q}\textbf{u}\|_{L^2}}{\|\textbf{u}\|_{B^{\sigma}_{2,1}}}.
\end{eqnarray*}
The third part $F^3_{q}$ is proceeded as follows:
\begin{eqnarray*}
F^3_{q}&=&-\triangle_{q}T_{\partial_{j}u^{j}}\tilde{m}\\
&=&-\sum_{|q-q'|\leqslant4}\triangle_{q}(S_{q'-1}\partial_{j}u^{j}\triangle_{q'}\tilde{m}),
\end{eqnarray*}
then
\begin{eqnarray*}
2^{q\sigma}\|F^3_{q}\|_{L^2}&\leqslant&  C\sum_{|q-q'|\leqslant4}2^{(q-q')\sigma}2^{q'\sigma}\|S_{q'-1}\partial_{j}u^{j}\triangle_{q'}\tilde{m}\|_{L^2}\\
&\leqslant& C\sum_{|q-q'|\leqslant4}2^{(q-q')\sigma}\|S_{q'-1}\partial_{j}u^{j}\|_{L^{\infty}}2^{q'(\sigma-1)}\|\triangle_{q'}\nabla\tilde{ m}\|_{L^2}\\
&\leqslant& Cc_{q3}\|\nabla
m\|_{B^{\sigma-1}_{2,1}}\|\textbf{u}\|_{B^{\sigma}_{2,1}},\ \ \
c_{q3}=\sum_{|q-q'|\leqslant4}\frac{2^{q'(\sigma-1)}\|\triangle_{q'}\nabla
m\|_{L^2}}{9\|\nabla m\|_{B^{\sigma-1}_{2,1}}}.
\end{eqnarray*}
(Here, we use the imbedding $B^{\sigma-1}_{2,1}\hookrightarrow
C_{0}$(continuous bounded functions which decay to zero at
infinity).)

 By the definition 2.2, we have
\begin{eqnarray*}F^4_{q}&=&\partial_{j}R(\tilde{m},\triangle_{q}u^{j})-\partial_{j}\triangle_{q}R(\tilde{m},u^{j})\\
&=&\sum_{|q-q'|\leqslant1}\partial_{j}(\triangle_{q'}\tilde{m}\tilde{\triangle}_{q'}\triangle_{q}u^{j})-\partial_{j}\triangle_{q}R(\tilde{m},u^{j})\\
&=&F^{4,1}_{q}+F^{4,2}_{q}.
\end{eqnarray*}
For the first term, using (4.8) only, we get
\begin{eqnarray*}
2^{q\sigma}\|F^{4,1}_{q}\|_{L^2}&\leqslant& C\|\nabla
m\|_{L^\infty}\sum_{|q-q'|\leqslant1}2^{(q-q')\sigma}2^{q'\sigma}\|\tilde{\triangle}_{q'}u^{j}\|_{L^2}\\
&\leqslant& Cc_{q4(1)}\|\nabla
m\|_{L^\infty}\|\textbf{u}\|_{B^{\sigma}_{2,1}},\ \ \
c_{4(1)}=\sum_{|q-q'|\leqslant1}\frac{2^{q'\sigma}\|\triangle_{q'}\textbf{u}\|_{L^2}}{4\|\textbf{u}\|_{B^{\sigma}_{2,1}}}.
\end{eqnarray*}
The second term is estimated as follows:
\begin{eqnarray*}
2^{q\sigma}\|F^{4,2}_{q}\|_{L^2}&=&2^{q\sigma}\|\partial_{j}\triangle_{q}R(\tilde{m},u^{j})\|_{L^2}\\
&\leqslant&C2^{q(\sigma+1)}\|\triangle_{q}R(\tilde{m},u^{j})\|_{L^2}\\
&\leqslant& Cc_{q4(2)}\|R(\tilde{m},u^{j})\|_{B^{\sigma+1}_{2,1}}\\
&\leqslant&
Cc_{q4(2)}\|m\|_{B^{\sigma}_{2,1}}\|\textbf{u}\|_{B^{\sigma}_{2,1}},\
\ \
c_{q4(2)}=\frac{2^{q(\sigma+1)}\|\triangle_{q}R(\tilde{m},u^{j})\|_{L^2}}{4\|R(\tilde{m},u^{j})\|_{B^{\sigma+1}_{2,1}}}.
\end{eqnarray*}
(Here, we use the result of continuity for the remainder, see
\cite{D} Proposition 1.4.2.\ \ \ $c_{q4}:=c_{q4(1)}+c_{q4(2)}$.)

\hspace{-5mm}For $F^5_{q}$, the same argument as  $F^4_{q}$, we can obtain
$$2^{q\sigma}\|F^{5}_{q}\|_{L^2}\leqslant
Cc_{q5}\|\nabla
m\|_{B^{\sigma-1}_{2,1}}\|\textbf{u}\|_{B^{\sigma}_{2,1}},\ \ \
c_{q5}=\Big(\sum_{|q-q'|\leqslant1}\frac{2^{q'\sigma}\|\triangle_{q'}\textbf{u}\|_{L^2}}{4\|\textbf{u}\|_{B^{\sigma}_{2,1}}}\Big)+\frac{2^{q\sigma}\|\triangle_{q}R(\partial_{j}\tilde{m},u^{j})\|_{L^2}}{4\|R(\partial_{j}\tilde{m},u^{j})\|_{B^{\sigma}_{2,1}}}.
$$
For
$F^6_{q}=\sum_{|q-q'|\leqslant1}[\triangle_{q}(\triangle_{-1}m\partial_{j}\triangle_{q'}u^{j})-\triangle_{-1}m\triangle_{q}\triangle_{q'}\partial_{j}u^{j}]
 \ (u^{j}=\sum_{q'}\triangle_{q'}u^{j}),$ applying first
order Taylor's formula, Young's inequality and (4.8), we have
\begin{eqnarray*}
2^{q\sigma}\|F^6_{q}\|_{L^2}&\leqslant&
C\sum_{|q-q'|\leqslant1}2^{(q-q')(\sigma-1)}\|\nabla
\triangle_{-1}m\|_{L^{\infty}}2^{q'\sigma}\|\triangle_{q'}\textbf{u}\|_{L^2}\\
&\leqslant& Cc_{q6}\|\nabla
m\|_{L^{\infty}}\|\textbf{u}\|_{B^{\sigma}_{2,1}},\ \ \ \
c_{q6}=\sum_{|q-q'|\leqslant1}\frac{2^{q'\sigma}\|\triangle_{q'}\textbf{u}\|_{L^2}}{3\|\textbf{u}\|_{B^{\sigma}_{2,1}}}.
\end{eqnarray*}
Adding above these inequalities together and choosing
$c_{q}=\frac{1}{6}\sum_{i=1}^6c_{qi}$,  we prove the
estimate (4.4).   \ \ \ \ \ $\square$ \\

 Now, we give the local
existence  result of solutions to (3.1)-(3.2).
\begin{prop}
Suppose that
$(m_{0},\textbf{u}_{0})\in{B^{\sigma}_{2,1}}$, then there exist
a time $T_{0}>0$ and a unique solution $(m, \textbf{u})$
of (3.1)-(3.2) such that $(m, \textbf{u})\in
C^{1}([0,T_{0}]\times \mathbb{R}^{N})$ and
$(m, \textbf{u})\in C([0,T_{0}],B^{\sigma}_{2,1})\cap
C^1([0,T_{0}],B^{\sigma-1}_{2,1})$.
\end{prop}
\hspace{-5mm}\textbf{Proof.} (\textbf{Existence})

 Let $U_{0}=(m_{0},
\textbf{u}_{0})^{T}\in{B^{\sigma}_{2,1}}$. There exists a sequence
$\{{U}^{k}_{0}\}:=\{(m_{0}^{k}, \textbf{u}^{k}_{0})^{T}\} \in H^{s}
(s>\sigma, \ s\in\mathbb{N})$ converging to $U_{0}$ in
$B^{\sigma}_{2,1}$ satisfying $
\|U^{k}_{0}\|_{B^{\sigma}_{2,1}}\leqslant
\|U_{0}\|_{B^{\sigma}_{2,1}}+1.$ We define a sequence
$\{U^{k}\}=\{(m^{k}, \textbf{u}^{k})^{T}\}$ solves the following
equations:
$$\cases{ m^{k}_{t}+\bar{\psi}\mbox{div}\textbf{u}^{k}=-\textbf{u}^{k}\cdot\nabla m^{k}-\frac{\gamma-1}{2}m^{k}\mbox{div}\textbf{u}^{k}\cr
 \textbf{u}^{k}_{t}+\bar{\psi}\nabla m^{k}+a\textbf{u}^{k}=-\textbf{u}^{k}\cdot\nabla\textbf{u}^{k}-\frac{\gamma-1}{2}m^{k}\nabla m^{k}}\eqno(4.9)$$
with the initial data
$$(
m^{k},\textbf{u}^{k})|_{t=0}=\Big(\frac{2}{\gamma-1}(\psi(n^{k}_{0})-\bar{\psi}),\textbf{u}^{k}_{0}\Big).\eqno(4.10)$$
It is easy to see (4.9) is a symmetric hyperbolic system on
$G=\{U^{k}:-\infty<\frac{\gamma-1}{2}m^{k}+\bar{\psi}<\infty\}$,
using Kato's classical results in \cite{K} or \cite{M}, we can get
the following local existence result: there exist a time $T_{k}>0$
and a solution $U^{k}$ of (4.9)-(4.10) such that
$$U^{k}\in C^{1}([0,T_{k}]\times \mathbb{R}^{N})$$ and $$U^{k}\in C([0,T_{k}],H^{s})\cap
C^{1}([0,T_{k}],H^{s-1}).$$ We define $[0,T^{*}_{k})$ is the
maximal interval of local existence for above solutions of
(4.9)-(4.10). According to the discussion in \cite{M}, we have the
blow-up criterion:
\begin{eqnarray*}T^{*}_{k}<\infty &\Leftrightarrow& \limsup_{t\rightarrow T^{*}_{k}}(\|U_{t}^{k}\|_{L^{\infty}}+\|\nabla
{U}^{k}\|_{L^{\infty}})=+\infty\mbox{\ \ \ or}\\
&&\mbox{for any compact subset} \ \ K \subset\subset G, \ \
U^{k}(x,t)\ \ \mbox{escapes}\ \  K \ \ \mbox{as} \ \ t\rightarrow
T^{*}_{k}.
\end{eqnarray*}
\textbf{Claim}: For $t\in [0,\min\{T^{*}_{k},T_{0}\}),$ we have
$\|U^{k}(t)\|_{B^{\sigma}_{2,1}}\leqslant 4\lambda_{0},$ where
$\lambda_{0}=\|U_{0}\|_{B^{\sigma}_{2,1}}+1$,
$T_{0}=1/(2\tilde{C}\lambda_{0})$ and $\tilde{C}$ is a positive
constant (independent of $k$) given in (4.12).

In fact, we have known
$$\|U^{k}_{0}\|_{B^{\sigma}_{2,1}}\leqslant
\lambda_{0}.$$  There exists a small $\bar{T}\in
(0,\min\{T^{*}_{k},T_{0}\})$ (owing to $U^{k}(t)\in
C([0,T_{k}],H^{s})$) such that
$$\sup_{t\in[0,\bar{T}]}\|U^{k}(t)\|_{B^{\sigma}_{2,1}}\leqslant 4\lambda_{0}.\eqno(4.11)$$
We can assume that (4.11) holds on arbitrary interval
$[0,T']\subset[0,\min\{T^{*}_{k},T_{0}\})$, then we shall show
(4.19) holds.\\Here, we don't need to consider the effect of
damping term. Therefore, by (4.2)-(4.4), H\"{o}lder's inequality
and Lemma 4.1, we can obtain
\begin{eqnarray*}\frac{d}{dt}(2^{q\sigma}\|\triangle_{q}U^{k}\|_{L^2})&\leqslant&
\tilde{C}\{\|\nabla\textbf{u}^{k}\|_{L^{\infty}}2^{q\sigma}(\|\triangle_{q}m^{k}\|_{L^2}+\|\triangle_{q}\textbf{u}^{k}\|_{L^2})+c_{q}\|\textbf{u}^{k}\|_{B^{\sigma}_{2,1}}\|m^{k}\|_{B^{\sigma}_{2,1}}
\\&&+c_{q}\|\nabla
\textbf{u}^{k}\|_{L^{\infty}}\|\textbf{u}^{k}\|_{B^{\sigma}_{2,1}}+\|\nabla
m^{k}\|_{L^{\infty}}2^{q\sigma}\|\triangle_{q}\textbf{u}^{k}\|_{L^2}\\&&+c_{q}\|\nabla
m^{k}\|_{L^{\infty}}\|m^{k}\|_{B^{\sigma}_{2,1}}+c_{q}\|m^{k}\|_{B^{\sigma}_{2,1}}\|\textbf{u}^{k}\|_{B^{\sigma}_{2,1}}\},\
\ t\in[0,T'].
\end{eqnarray*}$$\eqno(4.12)$$
$$(\|\triangle_{q}U^{k}\|^2_{L^2}:=\|\triangle_{q}m^{k}\|^2_{L^2}+\|\triangle_{q}\textbf{u}^{k}\|^2_{L^2},\,\,\|\triangle_{q}U^{k}\|_{L^2}\thickapprox\|\triangle_{q}m^{k}\|_{L^2}+\|\triangle_{q}\textbf{u}^{k}\|_{L^2}.)$$
Taking (4.12) $l^1$- norm, we obtain the
a-priori estimate of $U^{k}$:
$$\frac{d}{dt}\|U^{k}(t)\|_{B^{\sigma}_{2,1}}\leqslant
\tilde{C}\|U^{k}(t)\|^2_{B^{\sigma}_{2,1}},\ \ \ t\in[0,T'].
\eqno(4.13)$$ Integrating (4.13) on variable $t$ to get
\begin{eqnarray*}\|U^{k}(t)\|_{B^{\sigma}_{2,1}}&\leqslant&\|U^{k}_{0}\|_{B^{\sigma}_{2,1}}+\tilde{C}\int^{t}_{0}\|U^{k}(\varsigma)\|_{B^{\sigma}_{2,1}}^2d\varsigma\\
&\leqslant&(\|U_{0}\|_{B^{\sigma}_{2,1}}+1)+\tilde{C}\int^{t}_{0}\|U^{k}(\varsigma)\|_{B^{\sigma}_{2,1}}^2d\varsigma,
\ \ \ t\in[0,T'].\end{eqnarray*}$$\eqno(4.14)$$ Furthermore, we
have
$$\sup_{0\leqslant\varsigma\leqslant
t}\|U^{k}(\varsigma)\|_{B^{\sigma}_{2,1}}\leqslant
(\|U_{0}\|_{B^{\sigma}_{2,1}}+1)+\tilde{C}\int^{t}_{0}\sup_{0\leqslant\varsigma'\leqslant
\varsigma}\|U^{k}(\varsigma')\|_{B^{\sigma}_{2,1}}^2d\varsigma,\ \ \
t\in[0,T'].\eqno(4.15)$$ Set
$$\lambda_{1}(t)\equiv (\|U_{0}\|_{B^{\sigma}_{2,1}}+1)+\tilde{C}\int^{t}_{0}\sup_{0\leqslant\varsigma'\leqslant
\varsigma}\|U_{k}(\varsigma')\|_{B^{\sigma}_{2,1}}^2d\varsigma.$$
Then, we have
$$\frac{d}{dt}\lambda_{1}\leqslant \tilde{C}\lambda_{1}^2,\ \ \ \ \
\lambda_{1}(0)=\|U_{0}\|_{B^{\sigma}_{2,1}}+1,\ \
t\in[0,T'].\eqno(4.16)$$ Let $\lambda(t)$ solves Riccati equation:
$$\frac{d}{dt}\lambda=\tilde{C}\lambda^2,\ \ \ \ \
\lambda(0)=\|U_{0}\|_{B^{\sigma}_{2,1}}+1.\eqno(4.17)$$ The time
$T_{0}=1/(2\tilde{C}\lambda_{0})$ is less than the blow-up time
for (4.17). Then by solving the differential inequality (4.16), we
know $\lambda_{1}(t)\leqslant\lambda(t)$ for $t\in [0,T_{0}].$
Solving (4.17) yields
$$\lambda_{1}(t)\leqslant\frac{\|U_{0}\|_{B^{\sigma}_{2,1}}+1}{1-\tilde{C}t(\|U_{0}\|_{B^{\sigma}_{2,1}}+1)}=\lambda(t),
\ \ \ t\in[0,T_{0}]\ .\eqno(4.18)$$ Therefore, we see
$$
\sup_{
t\in[0,T']}\|U^{k}(t)\|_{B^{\sigma}_{2,1}}\leqslant2\lambda_{0}, \ \
\ [0,T']\subset[0,\min\{T^{*}_{k},T_{0}\}).\eqno(4.19)$$ Combining
with (4.11) and (4.19), by the continuum principle, we prove the
claim immediately.

 Furthermore, using Eqs.(4.9), we can conclude
$$\|U^{k}_{t}(t)\|_{B^{\sigma-1}_{2,1}}\leqslant \lambda_{0}',\ \ \ t\in[0,\min\{T^{*}_{k},T_{0}\}),$$
where $\lambda_{0}'$ is a positive constant only depending on the
initial data $U_{0}$. The blow-up criterion implies
$0<T_{0}<T^{*}_{k}$, so we have $0<T_{0}\leqslant
\inf_{k}T^{*}_{k}$.

That is, we find a positive time $T_{0}$ (only depending on the
initial data $U_{0}$) such that the approximative solution
sequence $\{U^{k}\}$ of (4.9)-(4.10) is uniform bounded in
$C([0,T_{0}],B^{\sigma}_{2,1})\cap
C^1([0,T_{0}],B^{\sigma-1}_{2,1}).$ Moreover, it
weak$^{*}$-converges (up to a subsequence) to some $U$ in
$L^{\infty}([0,T_{0}],B^{\sigma}_{2,1}),$ in terms of the
Banach-Alaoglu Theorem (see \cite{T} Remark 2 on p.180 in
\textsc{Triebel}, 1983). Because $\{U^{k}_{t}\} $ is also uniform
bounded in $C([0,T_{0}],B^{\sigma-1}_{2,1} )$(it
weak$^{*}$-converges to $U_{t}$ in
$L^{\infty}([0,T_{0}],B^{\sigma-1}_{2,1})$), then $\{U^{k}\}$ is
uniform bounded in $ \mathrm{Lip}([0,T_{0}],B^{\sigma-1}_{2,1})$,
hence uniform equicontinuous on $[0,T_{0}]$ valued in $
B^{\sigma-1}_{2,1}$. By Proposition 2.2, Ascoli-Arzela theorem and
Cantor diagonal process, we deduce that
$$\phi U^{k}\rightarrow \phi U\ \ \ \mbox{in}\ \ C([0,T_{0}],B^{\sigma-1}_{2,1}
)\ \ \mbox{as} \ \ k\rightarrow\infty, \ \ \ \mbox{for any}\ \
\phi\in\
 C_{c}^{\infty}.$$
The properties of strong convergence enable us to pass to the
limit in (4.9)-(4.10). Indeed, $U$ is a solution of (3.1)-(3.2).
Now, what remains is to check $U$ has the required regularity.
First, we already have known $U\in C([0,T_{0}],B^{\sigma-1}_{2,1}
)$, an interpolation argument insures $U\in
C([0,T_{0}],B^{\sigma'}_{2,1} )$ for any $\sigma'<\sigma.$
Furthermore, for any $q\in\mathbb{ N},\ S_{q}U\in
C([0,T_{0}],B^{\sigma}_{2,1}).$ Combining with (4.12)(throw off
the superscript $k$), we derive that $\{S_{q}U\}$ converges
uniformly to $U$ on $[0,T_{0}]$ valued in $B^{\sigma}_{2,1}$. This
achieves to prove that $U\in C([0,T_{0}],B^{\sigma}_{2,1}).$
Moreover, using Eqs.(3.1), we see that $U_{t}\in
C([0,T_{0}],B^{\sigma-1}_{2,1}),$ so $U(t,x)\in
C^{1}([0,T_{0}]\times \mathbb{R}^{N})$.\\

(\textbf{Uniqueness})

 Let $\tilde{m}=m_{1}-m_{2},\
\tilde{\textbf{u}}=\textbf{u}_{1}-\textbf{u}_{2},$ where
$U_{1}=(m_{1},\textbf{u}_{1})^{T},
U_{2}=(m_{2},\textbf{u}_{2})^{T}$ are two solutions for the system
(3.1)-(3.2) with the same initial data, respectively. Then
$\tilde{U}=(\tilde{m}, \tilde{\textbf{u}})$ satisfies the
following equations:
$$\cases{\tilde{m}_{t}+\bar{\psi}\mbox{div}\tilde{\textbf{u}}=-\textbf{u}_{1}\cdot\nabla \tilde{m}-\tilde{\textbf{u}}\nabla m_{2}-\frac{\gamma-1}{2}m_{1}\mbox{div}\tilde{\textbf{u}}-\frac{\gamma-1}{2}\tilde{m}\mbox{div}\textbf{u}_{2},\cr
 \tilde{\textbf{u}}_{t}+\bar{\psi}\nabla \tilde{m}+a\tilde{\textbf{u}}=-\textbf{u}_{1}\cdot\nabla\tilde{\textbf{u}}-\tilde{\textbf{u}}\nabla\textbf{u}_{2}-\frac{\gamma-1}{2}m_{1}\nabla\tilde{ m}-\frac{\gamma-1}{2}\tilde{ m}\nabla m_{2}. }\eqno(4.20)$$
Similar to the derivation of (4.13), we obtain the
following estimate:
$$ \|\tilde{U}(t)\|_{B^{\sigma-1}_{2,1}}\leqslant C
\int^{t}_{0}\|\tilde{U}(\varsigma)\|_{B^{\sigma-1}_{2,1}}(\|U_{1}(\varsigma)\|_{B^{\sigma}_{2,1}}+\|U_{2}(\varsigma)\|_{B^{\sigma}_{2,1}})d\varsigma\
, \ \ \ \mbox{for} \ \ t\in [0,T_{0}] .\eqno(4.21)$$ By Gronwall's
inequality, we conclude $\tilde{U}\equiv 0$.\ \ \ $\square$
\begin{rem}
If we add the assumption $\frac{\gamma-1}{2}m_{0}+\bar{\psi}>0$ in Proposition 4.1, defining the flow map $ X(t;s,x) $ of $\textbf{u}$ starting from $x\in \mathbb{R}^{N}$ at time $s\in[0,T_{0}]$ by
$$\frac{dX}{dt}=\textbf{u}(t,X(t;s,x)),\ \ \ X(t;s,x)|_{t=s}=x$$
and using the first equation of Eqs.(3.1), we can get
$$\frac{\gamma-1}{2}m(s,x)+\bar{\psi}>0 \ \ \mbox{for} \ (s,x)\in[0,T_{0}]\times \mathbb{R}^{N}. $$
\end{rem}

\section{ A-priori estimates and global existence }
In this section, we first give the proposition on the global
existence of classical solutions to (3.1)-(3.2).
\begin{prop} $(N\geqslant3)$ Suppose that $ U_{0}\in B^{\sigma+\varepsilon}_{2,2}$ ($\varepsilon>0$). There exists a positive constant $\delta_{2}<\frac{1}{2}\delta_{1}$ depending only on $A, \gamma, a$ and
$\bar{n}$ such that if
$$\|U(\cdot,0)\|^2_{B^{\sigma+\varepsilon}_{2,2}}+\|U_{t}(\cdot,0)\|^2_{B^{\sigma-1+\varepsilon}_{2,2}}\leqslant
\delta_{2}, $$ then there exists a unique global solution $U$ of
(3.1)-(3.2) satisfying
$$U\in C([0,\infty),B^{\sigma+\varepsilon}_{2,2})\cap C^1([0,\infty),B^{\sigma-1+\varepsilon}_{2,2})$$
and
\begin{eqnarray}&&\|U(\cdot,t)\|^2_{B^{\sigma+\varepsilon}_{2,2}}+\|U_{t}(\cdot,t)\|^2_{B^{\sigma-1+\varepsilon}_{2,2}}\nonumber\\&&+
\mu_{1}\int^{t}_{0}\Big(\|\textbf{u}(\cdot,\tau)\|^2_{B^{\sigma+\varepsilon}_{2,2}}+\|\nabla
m(\cdot,\tau)\|^2_{B^{\sigma-1+\varepsilon}_{2,2}}+\|U_{t}(\cdot,\tau)\|^2_{B^{\sigma-1+\varepsilon}_{2,2}}\Big)d\tau\nonumber\\&\leqslant&
\|U(\cdot,0)\|^2_{B^{\sigma+\varepsilon}_{2,2}}+\|U_{t}(\cdot,0)\|^2_{B^{\sigma-1+\varepsilon}_{2,2}},\
\ \ t\geqslant 0,\label{1}\end{eqnarray}  where $\delta_{1}$ and $
\mu_{1}$ are some positive constants given by Proposition 5.2,\
$U=(m,\textbf{u})$ and $U_{t}=(m_{t},\textbf{u}_{t})$.
\end{prop}
\begin{rem}
The energy estimate (5.1) implies the full solution $(m,\textbf{u})$
does not decay to zero in time exponentially.
\end{rem}
The proof of above proposition mainly depends on a crucial a-priori
estimate (Proposition 5.2). To do this, we need the following
lemmas.
\begin{lem} If $(m,\textbf{u})\in
C([0,T],B^{\sigma+\varepsilon}_{2,2})\cap
C^1([0,T],B^{\sigma-1+\varepsilon}_{2,2})$ is a solution of
Eqs.(3.1) for any given $T>0$, then
\begin{eqnarray}&&\frac{d}{dt}\Big(\|\triangle_{q}m_{t}\|^2_{L^2}+\|\triangle_{q}\textbf{u}_{t}\|^2_{L^2}\Big)+2a
\|\triangle_{q}\textbf{u}_{t}\|^2_{L^2}\nonumber\\&\leqslant&
2\|\textbf{u}_{t}\|_{L^{\infty}}(\|\triangle_{q}\nabla
m\|_{L^2}\|\triangle_{q}m_{t}\|_{L^2}+\|\triangle_{q}\nabla
\textbf{u}\|_{L^2}\|\triangle_{q}\textbf{u}_{t}\|_{L^2})\nonumber\\
&&+\|\nabla\textbf{u}\|_{L^{\infty}}(\|\triangle_{q}m_{t}\|^2_{L^2}+\|\triangle_{q}\textbf{u}_{t}\|^2_{L^2})+(\gamma-1)\|\nabla
m\|_{L^{\infty}}\|\triangle_{q}m_{t}\|_{L^2}\|\triangle_{q}\textbf{u}_{t}\|_{L^2}\nonumber\\
&&+(\gamma-1)\|m_{t}\|_{L^{\infty}}(\|\triangle_{q}\mathrm{div}\textbf{u}\|_{L^2}
\|\triangle_{q}m_{t}\|_{L^2}+\|\triangle_{q}\nabla
m\|_{L^2}\|\triangle_{q}\textbf{u}_{t}\|_{L^2})
+(\|[\textbf{u}_{t}, \triangle_{q}]\cdot\nabla
m\|_{L^2}\nonumber\\&&+\|[\textbf{u},\triangle_{q}]\cdot\nabla
m_{t}\|_{L^2})\|\triangle_{q}m_{t}\|_{L^2} +(\|[\textbf{u}_{t},
\triangle_{q}]\cdot\nabla \textbf{u}\|_{L^2}+\|[\textbf{u},
\triangle_{q}]\cdot\nabla
\textbf{u}_{t}\|_{L^2})\|\triangle_{q}\textbf{u}_{t}\|_{L^2}\nonumber\\
&&+(\gamma-1)(\|[m_{t},\triangle_{q}]\mathrm{div}\textbf{u}\|_{L^2}+\|[m,\triangle_{q}]\mathrm{div}\textbf{u}_{t}\|_{L^2})\|\triangle_{q}m_{t}\|_{L^2}\nonumber\\
&&+(\gamma-1)(\|[m_{t},\triangle_{q}]\nabla
m\|_{L^2}+\|[m,\triangle_{q}]\nabla
m_{t}\|_{L^2})\|\triangle_{q}\textbf{u}_{t}\|_{L^2}.   \label{2}
\end{eqnarray}
\end{lem}
\hspace{-5mm}\textbf{Proof.} By differentiating the first two
equations of Eqs.(3.1) with respect to the variable $t$ once,
integrating them over $\mathbb{R}^{N}$ after multiplying
$\triangle_{q}m_{t},\triangle_{q}\textbf{u}_{t}$, respectively,
similar to the derivation of (4.12), we can obtain (5.2) directly.
\ $\square$
\begin{lem}
The following estimates hold for any $m,\textbf{u}\in
C([0,T],B^{\sigma+\varepsilon}_{2,2})\cap C^1([0,T],B^{\sigma-1+\varepsilon}_{2,2}) (T>0)$:
\begin{eqnarray}2^{q(\sigma-1+\varepsilon)}\|[\textbf{u}_{t}, \triangle_{q}]\cdot\nabla
m\|_{L^2}\leqslant
Cc_{q}\|\textbf{u}_{t}\|_{B^{\sigma-1+\varepsilon}_{2,2}}\|\nabla
m\|_{B^{\sigma-1+\varepsilon}_{2,2}}, \label{3}\end{eqnarray}
\begin{eqnarray}2^{q(\sigma-1+\varepsilon)}\|[\textbf{u},\triangle_{q}]\cdot\nabla
m_{t}\|_{L^2}\leqslant
Cc_{q}\|\textbf{u}\|_{B^{\sigma+\varepsilon}_{2,2}}\|m_{t}\|_{B^{\sigma-1+\varepsilon}_{2,2}},\label{4}\end{eqnarray}
\begin{eqnarray}2^{q(\sigma-1+\varepsilon)}\|[\textbf{u}_{t}, \triangle_{q}]\cdot\nabla
\textbf{u}\|_{L^2}\leqslant C
c_{q}\|\textbf{u}_{t}\|_{B^{\sigma-1+\varepsilon}_{2,2}}\|\textbf{u}\|_{B^{\sigma+\varepsilon}_{2,2}},\label{5}\end{eqnarray}
\begin{eqnarray}2^{q(\sigma-1+\varepsilon)}\|[\textbf{u}, \triangle_{q}]\cdot\nabla
\textbf{u}_{t}\|_{L^2}\leqslant C
c_{q}\|\textbf{u}\|_{B^{\sigma+\varepsilon}_{2,2}}\|\textbf{u}_{t}\|_{B^{\sigma-1+\varepsilon}_{2,2}},
\label{6}\end{eqnarray}
\begin{eqnarray}2^{q(\sigma-1+\varepsilon)}\|[m_{t},\triangle_{q}]\mathrm{div}\textbf{u}\|_{L^2}\leqslant C
c_{q}\|m_{t}\|_{B^{\sigma-1+\varepsilon}_{2,2}}\|\textbf{u}\|_{B^{\sigma+\varepsilon}_{2,2}},\label{7}\end{eqnarray}
\begin{eqnarray}2^{q(\sigma-1+\varepsilon)}\|[m,\triangle_{q}]\mathrm{div}\textbf{u}_{t}\|_{L^2}\leqslant C
c_{q}\|m\|_{B^{\sigma+\varepsilon}_{2,2}}\|\textbf{u}_{t}\|_{B^{\sigma-1+\varepsilon}_{2,2}},\label{8}\end{eqnarray}
\begin{eqnarray}2^{q(\sigma-1+\varepsilon)}\|[m_{t},\triangle_{q}]\nabla m\|_{L^2}\leqslant C
c_{q}\|m_{t}\|_{B^{\sigma-1+\varepsilon}_{2,2}}\|\nabla
m\|_{B^{\sigma-1+\varepsilon}_{2,2}},\label{9}\end{eqnarray}
\begin{eqnarray}2^{q(\sigma-1+\varepsilon)}\|[m,\triangle_{q}]\nabla m_{t}\|_{L^2}\leqslant C
c_{q}\|m\|_{B^{\sigma+\varepsilon}_{2,2}}\|m_{t}\|_{B^{\sigma-1+\varepsilon}_{2,2}},\label{10}\end{eqnarray}
where $C$ denotes a harmless constant, $c_{q}$ denotes a sequence
such that $\|(c_{q})\|^2_{ {l^{2}}}\leqslant 1.$\end{lem}
\begin{rem} The proof is similar to that of Lemma 4.1, so we omit
 it here.
\end{rem}

In order to establish the differential inequality (5.30), we still
need some auxiliary estimates.

\begin{lem}If $(m,\textbf{u})\in C([0,T],B^{\sigma+\varepsilon}_{2,2})\cap
C^1([0,T],B^{\sigma-1+\varepsilon}_{2,2})$ is a solution of
Eqs.(3.1) for any given $T>0$, then
\begin{eqnarray}
\|\triangle_{q}m_{t}\|^2_{L^2}&\leqslant&
\Big(\bar{\psi}\bar{C}2^{q}\|\triangle_{q}\textbf{u}\|_{L^2}+\|\textbf{u}\|_{L^{\infty}}\|\triangle_{q}\nabla
m\|_{L^2}+\frac{\gamma-1}{2}\|m\|_{L^{\infty}}\|\triangle_{q}\mathrm{div}
\textbf{u}\|_{L^2}\nonumber\\&&+\|[\textbf{u},\triangle_{q}]\nabla
m\|_{L^2}+\frac{\gamma-1}{2}\|[m,\triangle_{q}]\mathrm{div}
\textbf{u}\|_{L^2}\Big)\|\triangle_{q}m_{t}\|_{L^2}, \label{11}
\end{eqnarray}
\begin{eqnarray}\bar{\psi}\|\triangle_{q}\nabla m\|^2_{L^2}&\leqslant& \Big(a\|\triangle_{q}\textbf{u}\|_{L^2}+\|\triangle_{q}\textbf{u}_{t}\|_{L^2}+\|\textbf{u}\|_{L^{\infty}}\|\triangle_{q}\nabla\textbf{u}\|_{L^2}+\|[\textbf{u},\triangle_{q}]\cdot\nabla \textbf{u}\|_{L^2}\nonumber\\&&+\frac{\gamma-1}{2}
\|m\|_{L^{\infty}}\|\triangle_{q}\nabla
m\|_{L^2}+\frac{\gamma-1}{2}\|[m,\triangle_{q}]\cdot\nabla
m\|_{L^2}\Big)\|\triangle_{q}\nabla m\|_{L^2},
\label{12}\end{eqnarray}
\end{lem} where the uniform constant $\bar{C}$ is independent of $A, \gamma,a$ and
$\bar{n}$.\newpage

\hspace{-5mm}\textbf{Proof.} (1) Using the first equation of
Eqs.(3.1), we have
\begin{eqnarray}m_{t}=-(\bar{\psi}\mbox{div}\textbf{u}+\textbf{u}\cdot\nabla
m+\frac{\gamma-1}{2}m\mbox{div}\textbf{u}).\label{13}\end{eqnarray}
By applying the operator $\triangle_{q}(q\geqslant-1)$ to (5.13),
integrating it over $\mathbb{R}^{N}$ after multiplying
$\triangle_{q}m_{t}$, we can get (5.11) only by H\"{o}lder's
inequality.

(2) Using the second equation of Eqs.(3.1), we get
\begin{eqnarray}\bar{\psi}\nabla
m=-(\textbf{u}_{t}+a\textbf{u}+\textbf{u}\cdot\nabla\textbf{u}+\frac{\gamma-1}{2}m\nabla
m).\label{14}\end{eqnarray} By applying the operator
$\triangle_{q}(q\geqslant-1)$ to (5.14), integrating it over
$\mathbb{R}^{N}$ after multiplying $\triangle_{q}\nabla m$, we can
get (5.12) immediately. \ \  $\square$

By Lemma 2.1, we have
 $$\|\triangle_{q}\nabla m\|_{L^2}\approx
 2^{q}\|\triangle_{q}m\|_{L^2}\ (q\geqslant0), $$ however, we can't get any estimates on  $\|\triangle_{-1}m\|_{L^2}$ according to  (5.12), furthermore, we can't obtain the total estimates
on $\|\triangle_{q} m\|_{L^2}(q\geqslant-1)$, which is the essential
difference with Euler-Poisson equation in \cite{FXZ}. That is why we
need the functional space $B^{\sigma+\varepsilon}_{2,2}$ to deal
with the global existence of classical solutions for (3.1)-(3.2).
Hence, we must modify those estimates in Lemma 4.1 and divide them
into the cases of the high and low frequency.
 \begin{lem}
    The following estimates hold for any $m,\textbf{u}\in B^{\sigma+\varepsilon}_{2,2}$:
\begin{eqnarray} 2^{q(\sigma+\varepsilon)}\|[\textbf{u},\triangle_{q}]\cdot\nabla \textbf{u}\|_{L^2}\leqslant Cc_{q}\|
   \textbf{u}\|_{B^{\sigma+\varepsilon}_{2,2}}\|\nabla\textbf{u}\|_{B^{\sigma-1+\varepsilon}_{2,2}}\ (q\geqslant -1),\label{15}\end{eqnarray}
 \begin{eqnarray}
  2^{q(\sigma+\varepsilon)}\|[m,\triangle_{q}]\nabla m\|_{L^2}\leqslant
 Cc_{q}\|\nabla m\|_{B^{\sigma+\varepsilon-1}_{2,2}}\|m\|_{B^{\sigma+\varepsilon}_{2,2}}\ (q\geqslant -1), \label{16}\end{eqnarray}
\begin{eqnarray}
 2^{q(\sigma+\varepsilon)}\|[\textbf{u},\triangle_{q}]\cdot\nabla m\|_{L^2}\leqslant
Cc_{q}\|\textbf{u}\|_{B^{\sigma+\varepsilon}_{2,2}}\|m\|_{B^{\sigma+\varepsilon}_{2,2}}\
(q\geqslant-1), \label{17} \end{eqnarray}
 \begin{eqnarray}
 2^{-(\sigma+\varepsilon)}\|[\textbf{u},\triangle_{-1}]\cdot\nabla m\|_{L^{\frac{2N}{N+2}}}\leqslant
Cc_{-1}\|\textbf{u}\|_{B^{\sigma+\varepsilon}_{2,2}}\|\nabla
m\|_{B^{\sigma-1+\varepsilon}_{2,2}} \ (N>2), \label{18}
 \end{eqnarray}
\begin{eqnarray}
   2^{q(\sigma+\varepsilon)}\|[m,\triangle_{q}]\mathrm{div}\textbf{u}\|_{L^2}\leqslant
   Cc_{q}\|m\|_{B^{\sigma+\varepsilon}_{2,2}}\|\textbf{u}\|_{B^{\sigma+\varepsilon}_{2,2}}\ (q\geqslant -1), \label{19}
\end{eqnarray}
   \begin{eqnarray}
   2^{-(\sigma+\varepsilon)}\|[m,\triangle_{-1}]\mathrm{div}\textbf{u}\|_{L^{\frac{2N}{N+2}}}\leqslant
   Cc_{-1}\|m\|_{B^{\sigma+\varepsilon}_{2,2}}\|\textbf{u}\|_{B^{\sigma+\varepsilon}_{2,2}} \ (N>2), \label{20}
  \end{eqnarray}
where $C$ denotes a  harmless  constant, $c_{q}(q\geqslant -1)$
denotes a sequence such that $\|(c_{q})\|^2_{ {l^{2}}}\leqslant
1.$
\end{lem}
 \begin{rem} The proof is similar to that of Lemma 4.1, so we also omit
 it here.
\end{rem}

Now, we give the crucial  a-priori estimate in the following
proposition.
\begin{prop}
There exist two positive constants $\delta_{1}$ and $\mu_{1}$
depending only on $A, \gamma, a$ and $\bar{n}$ such that for any
$T>0$, if  \begin{eqnarray} \sup_{0\leqslant t\leqslant T}\Big(\|U
(\cdot,t)\|^2_{B^{\sigma+\varepsilon}_{2,2}}+\|U_{t}(\cdot,t)\|^2_{B^{\sigma-1+\varepsilon}_{2,2}}\Big)\leqslant
\delta_{1},\label{21} \end{eqnarray} then
\begin{eqnarray}&&\|U(\cdot,t)\|^2_{B^{\sigma+\varepsilon}_{2,2}}+\|U_{t}(\cdot,t)\|^2_{B^{\sigma-1+\varepsilon}_{2,2}}\nonumber\\&&+
\mu_{1}\int^{t}_{0}\Big(\|\textbf{u}(\cdot,\tau)\|^2_{B^{\sigma+\varepsilon}_{2,2}}+\|\nabla
m(\cdot,\tau)\|^2_{B^{\sigma-1+\varepsilon}_{2,2}}+\|U_{t}(\cdot,\tau)\|^2_{B^{\sigma-1+\varepsilon}_{2,2}}\Big)d\tau\nonumber\\&\leqslant&
\|U(\cdot,0)\|^2_{B^{\sigma+\varepsilon}_{2,2}}+\|U_{t}(\cdot,0)\|^2_{B^{\sigma-1+\varepsilon}_{2,2}},\
\ \ t\geqslant 0.\label{1}\end{eqnarray}
\end{prop}
\hspace{-5mm}\textbf{Proof.} From the a-priori assumption (5.21), we
deduce
\begin{eqnarray}\sup_{0\leqslant t\leqslant T}\Big(\|U(\cdot,t)||_{L^{\infty}}+\|U(\cdot,t)||_{L^{N}}+\|\nabla
U(\cdot,t)||_{L^{\infty}}+\|\nabla
U(\cdot,t)||_{L^{N}}+\|U_{t}(\cdot,t)\|_{L^{\infty}}\Big)\leqslant
C\delta_{1}^{\frac{1}{2}}\ \ \ (N>2). \label{23}\end{eqnarray} In
addition, in order to obtain the global existence of $C^{1}$
solutions of the original system (1.1)-(1.3), we can choose
$0<\delta_{1}\leqslant\frac{\bar{\psi}^2}{(\gamma-1)^2C^2}$, then
$$\frac{\gamma-1}{2}m(t,x)+\bar{\psi}\geqslant
\frac{\bar{\psi}}{2}>0,\ \ \ (t,x)\in[0,T]\times\mathbb{R}^{N}.$$
From (4.2)-(4.3), by H\"{o}lder's inequality, we set
  \begin{eqnarray}
 I_{1,-1}&=&\|\textbf{u}\|_{L^{N}}\|\triangle_{-1}m\|_{L^{\frac{2N}{N-2}}}\|\triangle_{-1}\nabla m\|_{L^2}+\|\nabla\textbf{u}\|_{L^{\infty}}||\triangle_{-1}\textbf{u}\|^2_{L^2}
 +2\|[\textbf{u},\triangle_{-1}]\cdot\nabla
 m\|_{L^{\frac{2N}{N+2}}}\|\triangle_{-1}m\|_{L^{\frac{2N}{N-2}}}\nonumber\\&&+2\|[\textbf{u},\triangle_{-1}]\cdot\nabla
 \textbf{u}\|_{L^2}||\triangle_{-1}\textbf{u}\|_{L^2}+(\gamma-1)\|\nabla
 m\|_{L^{N}}||\triangle_{-1}m\|_{L^{\frac{2N}{N-2}}}\|\triangle_{-1}\textbf{u}\|_{L^2}\nonumber\\&&+(\gamma-1)\|[m,\triangle_{-1}]\nabla
 m\|_{L^2}||\triangle_{-1}\textbf{u}\|_{L^2}+(\gamma-1)\|[m,\triangle_{-1}]\mbox{div}\textbf{u}\|_{L^{\frac{2N}{N+2}}}||\triangle_{-1}m\|_{L^{\frac{2N}{N-2}}}\label{24}
 \end{eqnarray}
 and
\begin{eqnarray}
I_{1,q}&=&\|\nabla\textbf{u}\|_{L^{\infty}}(||\triangle_{q}m\|^2_{L^2}+||\triangle_{q}\textbf{u}\|^2_{L^2})
+2\|[\textbf{u},\triangle_{q}]\cdot\nabla
m\|_{L^2}||\triangle_{q}m\|_{L^2}+2\|[\textbf{u},\triangle_{q}]\cdot\nabla
\textbf{u}\|_{L^2}||\triangle_{q}\textbf{u}\|_{L^2}\nonumber\\&&+(\gamma-1)\|\nabla
m\|_{L^{\infty}}||\triangle_{q}m\|_{L^2}||\triangle_{q}\textbf{u}\|_{L^2}+(\gamma-1)\|[m,\triangle_{q}]\nabla
m\|_{L^2}||\triangle_{q}\textbf{u}\|_{L^2}\nonumber\\&&+(\gamma-1)\|[m,\triangle_{q}]\mbox{div}\textbf{u}\|_{L^2}||\triangle_{q}m\|_{L^2}\
 \ (q\geqslant0).\label{25}
\end{eqnarray}
$I_{2,q}$ denotes the right side of inequality $(5.1)\
(q\geqslant-1)$.

For the proof of Proposition 5.2, we are going to divided it into
the following two lemmas.
\begin{lem}$(q=-1)$
There exists a positive constant $\mu_{2}$ depending only on $A,
\gamma, a$ and $\bar{n}$ such that the following estimate holds:
\begin{eqnarray}
&&\frac{d}{dt}\Big(2^{-2(\sigma+\varepsilon)}\|\triangle_{-1}U\|^2_{L^2}+2^{-2(\sigma-1+\varepsilon)}\|\triangle_{-1}U_{t}\|^2_{L^2}\Big)\nonumber\\&&
+\mu_{2}\{2^{-2(\sigma+\varepsilon)}\|\triangle_{-1}\textbf{u}\|^2_{L^2}+2^{-2(\sigma-1+\varepsilon)}(\|\triangle_{-1}\nabla m\|^2_{L^2}+\|\triangle_{-1}U_{t}\|^2_{L^2})\}\nonumber\\
&\leqslant& C\{(\|\textbf{u}\|_{L^{N}}+\|\nabla
m\|_{L^{N}}+\|\nabla\textbf{u}\|_{L^{\infty}}+\|U\|_{B^{\sigma+\varepsilon}_{2,2}})(2^{-2(\sigma+\varepsilon-1)}\|\triangle_{-1}\nabla
m\|^2_{L^2}\nonumber\\&&+2^{-2(\sigma+\varepsilon)}\|\triangle_{-1}\textbf{u}\|^2_{L^2}
 +c^2_{-1}\|\textbf{u}\|^2_{B^{\sigma+\varepsilon}_{2,2}}+c^2_{-1}\|\nabla U\|^2_{B^{\sigma-1+\varepsilon}_{2,2}})\}+J_{-1},\label{26}\end{eqnarray}
where
\begin{eqnarray}
J_{q}&=&C(\|U\|_{L^{\infty}}+\|\nabla
U\|_{L^{\infty}}+\|U_{t}\|_{L^{\infty}})2^{2q(\sigma-1+\varepsilon)}(\|\triangle_{q}\nabla
U\|^2_{L^2}+\|\triangle_{q}U_{t}\|^2_{L^2})\nonumber\\&&+(\|U\|_{B^{\sigma+\varepsilon}_{2,2}}+\|U_{t}\|_{B^{\sigma-1+\varepsilon}_{2,2}})\{2^{2q(\sigma-1+\varepsilon)}(\|\triangle_{q}\nabla
m\|^2_{L^2}+\|\triangle_{q}U_{t}\|^2_{L^2})\nonumber\\&&+c_{q}^2(\|\textbf{u}\|^2_{B^{\sigma+\varepsilon}_{2,2}}+\|\nabla
U\|^2_{B^{\sigma-1+\varepsilon}_{2,2}}+\|U_{t}\|^2_{B^{\sigma-1+\varepsilon}_{2,2}})\}\
(q\geqslant-1).\nonumber
\end{eqnarray}
\end{lem}
\hspace{-5mm}\textbf{Proof of Lemma 5.5.} Combining (4.2)-(4.3),
Lemma 5.1 and 5.3, we have
\begin{eqnarray}
&&\frac{d}{dt}\Big\{2^{-2}(\|\triangle_{-1}m\|^2_{L^2}+\|\triangle_{-1}\textbf{u}\|^2_{L^2})+(\|\triangle_{-1}m_{t}\|^2_{L^2}+\|\triangle_{-1}\textbf{u}_{t}\|^2_{L^2})\Big\}
\nonumber\\&&+\beta_{1}\bar{\psi}\|\triangle_{-1}\nabla
m\|^2_{L^2}+2a2^{-2}\|\triangle_{-1}\textbf{u}\|^2_{L^2}
+\beta_{2}\|\triangle_{-1}m_{t}\|^2_{L^2}+2a\|\triangle_{-1}\textbf{u}_{t}\|^2_{L^2}
\nonumber\\&\leqslant& 2^{-2}I_{1,-1}
+I_{2,-1}+\beta_{1}\Big(a\|\triangle_{-1}\textbf{u}\|_{L^2}+\|\triangle_{-1}\textbf{u}_{t}\|_{L^2}+\|\textbf{u}\|_{L^{\infty}}\|\triangle_{-1}\nabla
\textbf{u}\|_{L^2}\nonumber\\&&+\|[\textbf{u},\triangle_{-1}]\nabla
\textbf{u}\|_{L^2}+\frac{\gamma-1}{2}\|m\|_{L^{\infty}}\|\triangle_{-1}\nabla
m\|_{L^2}+\frac{\gamma-1}{2}\|[m,\triangle_{-1}]\nabla
m\|_{L^2}\Big)\|\triangle_{-1}\nabla
m\|_{L^2}\nonumber\\&&+\beta_{2}
\Big(\bar{\psi}\bar{C}2^{-1}\|\triangle_{-1}\textbf{u}\|_{L^2}+\|\textbf{u}\|_{L^{\infty}}\|\triangle_{-1}\nabla
m\|_{L^2}+\frac{\gamma-1}{2}\|m\|_{L^{\infty}}\|\triangle_{-1}\mathrm{div}
\textbf{u}\|_{L^2}\nonumber\\&&+\frac{\gamma-1}{2}\|[m,\triangle_{-1}]\mathrm{div}
\textbf{u}\|_{L^2}+\|[\textbf{u},\triangle_{-1}]\nabla
m\|_{L^2}\Big)\|\triangle_{-1}m_{t}\|_{L^2},\label{27}
\end{eqnarray}
where two positive constants $\beta_{1},\beta_{2}$ satisfy
$$\beta_{1}\leqslant\min\Big\{\frac{\bar{\psi}}{8a}, \
\bar{\psi}a\Big\}\ \mbox{and}\, \
\beta_{2}\leqslant\frac{a}{\bar{\psi}^2\bar{C}^2},\ \
\mbox{respectively}.$$ We introduce them in order to eliminate
quadratic terms in the right side of (5.27). First, we notice that
there are no quadratic terms in $I_{1,-1}$ and $I_{2,-1}$. The first
quadratic term can be handled directly by Young's inequality:
\begin{eqnarray*}
&&\beta_{1}a\|\triangle_{-1}\textbf{u}\|_{L^2}\|\triangle_{-1}\nabla m\|_{L^2}\\
&=&\sqrt{\frac{\beta_{1}}{\bar{\psi}}}a\|\triangle_{-1}\textbf{u}\|_{L^2}\cdot
\sqrt{\beta_{1}\bar{\psi}}\|\triangle_{-1}\nabla m\|_{L^2}\\
&\leqslant&\frac{\beta_{1}a^2}{\bar{\psi}}\|\triangle_{-1}\textbf{u}\|^2_{L^2}+\frac{1}{4}\beta_{1}\bar{\psi}\|\triangle_{-1}\nabla m\|^2_{L^2}\\
&\leqslant&\frac{1}{8}a\|\triangle_{-1}\textbf{u}\|^2_{L^2}+\frac{1}{4}\beta_{1}\bar{\psi}\|\triangle_{-1}\nabla
m\|^2_{L^2}.\end{eqnarray*} The remainder quadratic terms in the
right side of (5.27) are estimated similarly as follows:
\begin{eqnarray*}\beta_{1}\|\triangle_{-1}\textbf{u}_{t}\|_{L^2}\|\triangle_{-1}\nabla m\|_{L^2}&\leqslant&\frac{\beta_{1}}{\bar{\psi}}\|\triangle_{-1}\textbf{u}_{t}\|^2_{L^2}+\frac{1}{4}\beta_{1}\bar{\psi}\|\triangle_{-1}\nabla m\|^2_{L^2}\\
&\leqslant&a\|\triangle_{-1}\textbf{u}_{t}\|^2_{L^2}+\frac{1}{4}\beta_{1}\bar{\psi}\|\triangle_{-1}\nabla m\|^2_{L^2};
\end{eqnarray*}
\begin{eqnarray*}\frac{1}{2}\beta_{2}\bar{\psi}\bar{C}\|\triangle_{-1}\textbf{u}\|_{L^2}\|\triangle_{-1}m_{t}\|_{L^2}
&\leqslant&\frac{1}{8}\beta_{2}\bar{\psi}^2\bar{C}^2\|\triangle_{-1}\textbf{u}\|^2_{L^2}+\frac{1}{2}\beta_{2}\|\triangle_{-1}m_{t}\|^2_{L^2},\\
&\leqslant&\frac{1}{8}a\|\triangle_{-1}\textbf{u}\|^2_{L^2}+\frac{1}{2}\beta_{2}\|\triangle_{-1}m_{t}\|^2_{L^2}.
\end{eqnarray*}
 Then (5.27) becomes into
 \begin{eqnarray}
 &&\frac{d}{dt}\Big\{2^{-2}(\|\triangle_{-1}m\|^2_{L^2}+\|\triangle_{-1}\textbf{u}\|^2_{L^2})+(\|\triangle_{-1}m_{t}\|^2_{L^2}+\|\triangle_{-1}\textbf{u}_{t}\|^2_{L^2})\Big\}
 \nonumber\\&&+\frac{1}{2}\beta_{1}\bar{\psi}\|\triangle_{-1}\nabla m\|^2_{L^2}+\frac{1}{4}a\|\triangle_{-1}\textbf{u}\|^2_{L^2}
 +\frac{1}{2}\beta_{2}\|\triangle_{-1}m_{t}\|^2_{L^2}+a\|\triangle_{-1}\textbf{u}_{t}\|^2_{L^2}
 \nonumber\\&\leqslant& 2^{-2}I_{1,-1}
 +I_{2,-1}+\beta_{1}\Big(\|\textbf{u}\|_{L^{\infty}}\|\triangle_{-1}\nabla
 \textbf{u}\|_{L^2}+\|[\textbf{u},\triangle_{-1}]\nabla
 \textbf{u}\|_{L^2}\nonumber\\&&+\frac{\gamma-1}{2}\|m\|_{L^{\infty}}\|\triangle_{-1}\nabla
 m\|_{L^2}+\frac{\gamma-1}{2}\|[m,\triangle_{-1}]\nabla
 m\|_{L^2}\Big)\|\triangle_{-1}\nabla m\|_{L^2}\nonumber\\&&+\beta_{2}
 \Big(\|\textbf{u}\|_{L^{\infty}}\|\triangle_{-1}\nabla
 m\|_{L^2}+\frac{\gamma-1}{2}\|m\|_{L^{\infty}}\|\triangle_{-1}\mathrm{div}
 \textbf{u}\|_{L^2}\nonumber\\&&+\frac{\gamma-1}{2}\|[m,\triangle_{-1}]\mathrm{div}
 \textbf{u}\|_{L^2}+\|[\textbf{u},\triangle_{-1}]\nabla
m\|_{L^2}\Big)\|\triangle_{-1}m_{t}\|_{L^2}.\label{28}
\end{eqnarray}
 Multiplying (5.28) by $2^{-2(\sigma-1+\varepsilon)}$ and combining Lemma 5.2, 5.4, we can get  (5.26) with the aid of
 Gagliardo-Nirenberg-Sobolev inequality $\Big(\|\triangle_{-1}m\|_{L^{\frac{2N}{N-2}}}\leqslant C
\|\triangle_{-1}\nabla m\|_{L^2}\ (N>2)\Big)$ and  Young's
inequality.\ \ \ \ $ \square$

For the case of high frequency ($q \geqslant 0$), we also have the
following a-priori estimate in a similar way:
\begin{lem}$(q\geqslant 0)$
There exists a positive constant $\mu_{3}$ depending only on $A,
\gamma, a$ and $\bar{n}$ such that the following estimate holds:
\begin{eqnarray}
&&\frac{d}{dt}\Big(2^{2q(\sigma+\varepsilon)}\|\triangle_{q}U\|^2_{L^2}+2^{2q(\sigma-1+\varepsilon)}\|\triangle_{q}U_{t}\|^2_{L^2}\Big)\nonumber\\&&
+\mu_{3}\{2^{2q(\sigma+\varepsilon)}\|\triangle_{q}\textbf{u}\|^2_{L^2}+2^{2q(\sigma-1+\varepsilon)}(\|\triangle_{q}\nabla m\|^2_{L^2}+\|\triangle_{q}U_{t}\|^2_{L^2})\}\nonumber\\
&\leqslant& C\{(\|\nabla
U\|_{L^{\infty}}+\|U\|_{B^{\sigma+\varepsilon}_{2,2}})(2^{2q(\sigma-1+\varepsilon)}\|\triangle_{q}\nabla
m\|^2_{L^2}+2^{2q(\sigma+\varepsilon)}\|\triangle_{q}\textbf{u}\|^2_{L^2}
 \nonumber\\&&+c_{q}^2\|\textbf{u}\|^2_{B^{\sigma+\varepsilon}_{2,2}}+c_{q}^2\|\nabla
U\|^2_{B^{\sigma-1+\varepsilon}_{2,2}})\} +J_{q}, \label{29}
\end{eqnarray}
where $J_{q}$ is defined in Lemma 5.5.
\end{lem}
  Summing (5.29) on $q \in \mathbb{N}\cup\{0\}$ and adding (5.26) together, according to a-priori assumption
(5.23), we get the following differential inequality:
\begin{eqnarray} &&\frac{d}{dt}\Big(\|U(\cdot,t)\|^2_{B^{\sigma+\varepsilon}_{2,2}}+\|U_{t}(\cdot,t)\|^2_{B^{\sigma-1+\varepsilon}_{2,2}}\Big)\nonumber\\&&+
\mu_{4}\Big(\|\textbf{u}(\cdot,t)\|^2_{B^{\sigma+\varepsilon}_{2,2}}+\|\nabla
m\|^2_{B^{\sigma-1+\varepsilon}_{2,2}}+\|U_{t}(\cdot,t)\|^2_{B^{\sigma-1+\varepsilon}_{2,2}}\Big)\nonumber\\&&\hspace{-4mm}\leqslant
C\delta_{1}^{\frac{1}{2}}\Big(\|\textbf{u}(\cdot,t)\|^2_{B^{\sigma+\varepsilon}_{2,2}}+\|\nabla
m\|^2_{B^{\sigma-1+\varepsilon}_{2,2}}+\|U_{t}(\cdot,t)\|^2_{B^{\sigma-1+\varepsilon}_{2,2}}\Big),\label{30}\
\end{eqnarray}
where the constant $\mu_{4}$ depends only on $A, \gamma, a$ and
$\bar{n}$. Furthermore, choosing
$\delta_{1}\leqslant\min\{\frac{\mu_{4}^2}{4C^2},\frac{\bar{\psi}^2}{(\gamma-1)^2C^2}\}$,
we conclude the proof
of Proposition 5.2 with $\mu_{1}=\frac{\mu_{4}}{2}$.\ \ \ $\square$\\

 \hspace{-5mm}\textbf{Proof of
Proposition 5.1.} In fact,\ Proposition 4.1 also holds on the
framework of functional space $B^{\sigma+\varepsilon}_{2,2}$. From
the assumption
$$\|U(\cdot,0)\|^2_{B^{\sigma+\varepsilon}_{2,2}}+\|U_{t}(\cdot,0)\|^2_{B^{\sigma-1+\varepsilon}_{2,2}}\leqslant
\delta_{2},$$ we can determine a time $T_{1}>0 \ (T_{1}<T_{0})$ such
that
$$\|U(\cdot,t)\|^2_{B^{\sigma+\varepsilon}_{2,2}}+\|U_{t}(\cdot,t)\|^2_{B^{\sigma-1+\varepsilon}_{2,2}}\leqslant2
\delta_{2}, \ \ \ \mbox{for all} \ \ t\in [0,T_{1}].$$
\textbf{Claim}:  One can choose a positive constant $\delta_{2}$
satisfying $\delta_{2}<\frac{1}{2}\delta_{1}$ such that
\begin{eqnarray}\|U(\cdot,t)\|^2_{B^{\sigma+\varepsilon}_{2,2}}+\|U_{t}(\cdot,t)\|^2_{B^{\sigma-1+\varepsilon}_{2,2}}<
\delta_{1}, \ \ \ \mbox{for all} \ \ t\in
[0,T_{0}].\label{31}\end{eqnarray} In fact, otherwise, we can assume
that there exists a time $T_{2}\
 (T_{1}<T_{2}\leqslant T_{0})$ such that (5.31) is satisfied for all
$t\in[0,T_{2})$ and
\begin{eqnarray}\|U(\cdot,T_{2})\|^2_{B^{\sigma+\varepsilon}_{2,2}}+\|U_{t}(\cdot,T_{2})\|^2_{B^{\sigma-1+\varepsilon}_{2,2}}=
\delta_{1},\label{32}\end{eqnarray} since we see (5.31) is satisfied
as $t\in[0,T_{1}]$ for such choice of $\delta_{2}$.\\
 By Proposition 5.2, for all
$t\in[0,T^{k}](T^{k}\rightarrow T_{2})$
$$\|U(\cdot,t)\|^2_{B^{\sigma+\varepsilon}_{2,2}}+\|U_{t}(\cdot,t)\|^2_{B^{\sigma-1+\varepsilon}_{2,2}}\leqslant
\|U(\cdot,0)\|^2_{B^{\sigma+\varepsilon}_{2,2}}+\|U_{t}(\cdot,0)\|^2_{B^{\sigma-1+\varepsilon}_{2,2}}.$$
In particular,
$$
\|U(\cdot,T^{k})\|^2_{B^{\sigma+\varepsilon}_{2,2}}+\|U_{t}(\cdot,T^{k})\|^2_{B^{\sigma-1+\varepsilon}_{2,2}}\leqslant
\|U(\cdot,0)\|^2_{B^{\sigma+\varepsilon}_{2,2}}+\|U_{t}(\cdot,0)\|^2_{B^{\sigma-1+\varepsilon}_{2,2}}.
$$
By the continuity on $t\in[0,T_{0}]$, we get
$$\|U(\cdot,T_{2})\|^2_{B^{\sigma+\varepsilon}_{2,2}}+\|U_{t}(\cdot,T_{2})\|^2_{B^{\sigma-1+\varepsilon}_{2,2}}\leqslant
\delta_{2}<\delta_{1},$$ which contradicts (5.32). Hence, (5.31) is
true.  From Proposition 4.1 and 5.2, we can prove Proposition 5.1 by
using the standard boot-strap argument. \ \ \ $\square$

By Besov imbedding property, $(m,\textbf{u})\in
C^{1}([0,\infty)\times \mathbb{R}^{N})$ solves (3.1)-(3.2). The
choice of $\delta_{1}$ is sufficient to ensure
$\frac{\gamma-1}{2}m+\bar{\psi}>0$. According to Remark 3.1, we
deduce that $(n,\textbf{u})\in C^{1}([0,\infty)\times
\mathbb{R}^{N})$ solves (1.1)-(1.3) with $n>0$. Furthermore, we
attain the main result (Theorem 1.2) in this paper.\\

In what follows, we state a direct consequence of Proposition 5.1.
\begin{cor}
Let $(m, \textbf{u})$ be the solution in Proposition 5.1, we have ($
\sigma=1+\frac{N}{2},\ \varepsilon'<\varepsilon.)$
$$\|m(\cdot,t)\|_{B^{\sigma-1+\varepsilon'}_{p,2}}\rightarrow 0 \ \ (p=\frac{2N}{N-2}),\ \ \ \ \|\textbf{u}(\cdot,t)\|_{B^{\sigma+\varepsilon'}_{2,2}}\rightarrow 0,\ \ \mbox{as}\ \ t\rightarrow +\infty.   $$
\end{cor}
\hspace{-5mm}\textbf{Proof.} Because of the similar argument, we
show the former only. From the energy estimate in Proposition 5.1,
we get
$$\nabla m\in L^2_{t}B^{\sigma-1+\varepsilon}_{2,2},\ \ \ \ \nabla m_{t}\in L^2_{t}B^{\sigma-2+\varepsilon}_{2,2}.$$
Set $$H(t)=\|\nabla
m(\cdot,t)\|^2_{B^{\sigma-2+\varepsilon}_{2,2}}\in L^1_{t}.$$ After
an easy computation, we have
$$ \frac{d}{dt}H(t)\leqslant \|\nabla
m(\cdot,t)\|^2_{B^{\sigma-2+\varepsilon}_{2,2}}+\|\nabla
m_{t}(\cdot,t)\|^2_{B^{\sigma-2+\varepsilon}_{2,2}}\in L^1_{t}.$$
Hence, $H(t)\rightarrow 0$ as $t\rightarrow +\infty$. Since $m(t,x)$
is bounded in $C([0,\infty),B^{\sigma+\varepsilon}_{2,2})$, by
interpolation argument, we can obtain $(\varepsilon'<\varepsilon )$
$$\|\nabla m(\cdot,t)\|_{B^{\sigma-1+\varepsilon'}_{2,2}}\rightarrow
0, \ \ \mbox{as}\ \ t\rightarrow +\infty,$$ which completes the
proof after using Gagliardo-Nirenberg-Sobolev
inequality.\ \ \ $\square$\\

Finally, we show the exponential decay of the vorticity.

\hspace{-5mm}\textbf{Proof of Theorem 1.3.} When $N=3$, the curl
of the velocity equation in Eqs.(1.1) gives
$$\partial_{t}\omega+a\omega+\textbf{u}\cdot\nabla \omega-\omega\cdot\nabla\textbf{u}=0.$$ Hence,
\begin{eqnarray}\frac{1}{2}\frac{d}{dt}\|\triangle_{q}\omega\|^2_{L^2}+a\|\triangle_{q}\omega\|^2_{L^2}\leqslant
C_{0}(\|\nabla\textbf{u}\|_{L^{\infty}}\|\triangle_{q}\omega\|_{L^2}+\|\omega\|_{L^{\infty}}\|\triangle_{q}\nabla\textbf{u}\|_{L^2}+c_{q}\|\nabla
\textbf{u}\|_{B^{\sigma-1}_{2,1}}\|\omega\|_{B^{\sigma-1}_{2,1}})\|\triangle_{q}\omega\|_{L^2}.
\hspace{1.5mm}\label{33}\end{eqnarray} Dividing (5.33) by
$\|\triangle_{q}\omega\|_{L^2}$ and summing it on
$q\geqslant-1(q\in \mathbb{Z})$  after multiplying the factor
$2^{q(\sigma-1)}$, from Theorem 1.2, we have
\begin{eqnarray*}&&\frac{1}{2}\frac{d}{dt}\|\omega(\cdot,t)\|_{B^{\sigma-1}_{2,1}}+a\|\omega(\cdot,t)\|_{B^{\sigma-1}_{2,1}}\\&\leqslant&
C_{0}\|\textbf{u}(\cdot,t)\|_{B^{\sigma+\varepsilon}_{2,2}}\|\omega(\cdot,t)\|_{B^{\sigma-1}_{2,1}}\\&\leqslant&
C_{0}\Big(\|(n-\bar{n},\textbf{u})(\cdot,0)\|^2_{B^{\sigma+\varepsilon}_{2,2}}+\|(n_{t},\textbf{u}_{t})(\cdot,0)\|^2_{B^{\sigma-1+\varepsilon}_{2,2}}\Big)^{\frac{1}{2}}\|\omega(\cdot,t)\|_{B^{\sigma-1}_{2,1}}
\\&\leqslant&C_{0}\min\{\delta_{0}^{\frac{1}{2}}, \frac{a}{2C_{0}}\}\|\omega(\cdot,t)\|_{B^{\sigma-1}_{2,1}}\\&\leqslant&\frac{1}{2}a\|\omega(\cdot,t)\|_{B^{\sigma-1}_{2,1}}.\end{eqnarray*}
Therefore, we obtain the exponential decay of
$\|\omega(\cdot,t)\|_{B^{\sigma-1}_{2,1}}$ with $\mu'_{0}=a$.\ \ \
$\square$

\end{document}